\documentclass[10pt,oneside,reqno]{amsart}
\usepackage[T1]{fontenc}
\usepackage[utf8]{inputenc}
\usepackage[british]{babel}
\usepackage{bm,blindtext,amsmath,amssymb,amsthm,mathtools,tikz-cd,mathrsfs,cancel,xcolor,wrapfig,comment,array,pbox,booktabs,slashed,graphicx,subfiles,multirow,multicol}
\usepackage{enumitem}
\usepackage{tabularray}
\graphicspath{ {./images/} }

\usepackage[colorlinks=true,allcolors=blue,hyperindex,breaklinks,hypertexnames=false]{hyperref}

\calclayout

\numberwithin{equation}{section}
\theoremstyle{plain}
\newtheorem{theorem}{Theorem}[section]
\newtheorem{lemma}[theorem]{Lemma}
\newtheorem{prop}[theorem]{Proposition}
\newtheorem{cor}[theorem]{Corollary}
\newtheorem{conj}{Conjecture}

\newcounter{step}[theorem]
\renewcommand{\thestep}{\arabic{step}}

\theoremstyle{definition}
\newtheorem{defn}[theorem]{Definition}
\newtheorem{ex}[theorem]{Example}

\newtheorem{remark}[theorem]{Remark}

\newcommand{\git}{\mathbin{
  \mathchoice{/\mkern-6mu/}
    {/\mkern-6mu/}
    {/\mkern-5mu/}
    {/\mkern-5mu/}}}

\renewcommand{\a}{\alpha}
\renewcommand{\b}{\beta}

\newcommand{\N}{\mathcal{N}}
\newcommand{\MM}{\mathscr{M}}
\newcommand{\W}{\mathcal{W}}
\newcommand{\Z}{\mathbb{Z}}

\newcommand{\A}{\mathcal{A}}

\newcommand{\M}{\mathcal{M}}

\newcommand{\DD}{\mathbb{D}}

\newcommand{\C}{\mathbb{C}}
\renewcommand{\P}{\mathbb{P}}

\newcommand{\ra}{\rightarrow}

\newcommand{\Spec}{\textnormal{Spec}}

\newcommand{\gr}{\textnormal{gr}}
\newcommand{\tr}{\textnormal{tr}}

\newcommand{\GL}{\textnormal{GL}}

\newcommand{\Quot}{\textnormal{Quot}}

\newcommand{\Bl}{\textnormal{Bl}}
\newcommand{\SL}{\textnormal{SL}}
\newcommand{\SU}{\textnormal{SU}}

\newcommand{\PGL}{\textnormal{PGL}}

\newcommand{\e}{\varepsilon}

\newcommand{\id}{\textnormal{id}}

\newcommand{\codim}{\textnormal{codim}}
\newcommand{\im}{\textnormal{im}}
\newcommand{\End}{\textnormal{End}}

\newcommand{\Hom}{\textnormal{Hom}}

\newcommand{\Ext}{\textnormal{Ext}}

\newcommand{\Pic}{\textnormal{Pic}}

\newcommand{\ev}{\textnormal{ev}}

\newcommand{\Fix}{\textnormal{Fix}}
\newcommand{\Aut}{\textnormal{Aut}}

\newcommand{\Spin}{\textnormal{Spin}}

\newcommand{\pr}{\textnormal{pr}}
\newcommand{\E}{\mathcal{E}}

\newcommand{\F}{\mathcal{F}}

\renewcommand{\O}{\mathcal{O}}

\newcommand{\del}{\partial}
\newcommand{\delbar}{\overline{\partial}}

\title{Multi-monopoles and monodromy around the wobbly divisor}

\author{Brad Wilson}
\address{Department of Mathematics, University College London, London, WC1H 0AY}
\email{brad.wilson.23@ucl.ac.uk}

\begin{document}

\begin{abstract}
    The multi-monopole equations are a generalisation of the Seiberg--Witten equations where the spinor is coupled to an auxiliary vector bundle. The signed count of multi-monopoles on a 3-manifold can vary according to a new type of wall-crossing phenomenon. For products $S^1\times C$, where $C$ is a closed Riemann surface, we relate the wall and chamber structure to loops in the complement of the wobbly divisor of $C$. Such a loop induces a monodromy map on the multi-vortex moduli space, and the count of its fixed points is expected to agree with the count of multi-monopoles on $S^1\times C$ in the given chamber. By studying singularities of framed multi-vortex moduli spaces, we are able to understand these monodromy maps for various loops. Counting fixed points suggests that product 3-manifolds display incredibly rich wall-crossing behaviour. 
\end{abstract}

\maketitle

\section{Introduction}
\label{intro}
In this article, we take an algebro-geometric approach towards the study of wall-crossing phenomena arising in 2-spinor Seiberg--Witten theory on 3-manifolds. Before stating our main results, we review the idea of wall-crossing and how the author's previous work \cite{wilson} relates it to questions about the monodromy of certain families of moduli spaces. 

\subsection{Multi-monopoles on 3-manifolds}\label{intro1}

The multi-monopole equations are a generalisation of the Seiberg--Witten equations on a closed 3-manifold $Y$, where the spinor field is coupled to an auxiliary rank $N$ vector bundle. In this paper, we only consider the $N=2$ case, which already captures the main features of the theory. The equations depend on a spin$^c$ structure $\mathfrak{s}\in\Spin^c(Y)$ and a pair $\mathfrak{p}=(h,B)$, where $h$ is a Riemannian metric on $Y$ and $B$ is a $\SU(2)$-connection on the auxiliary bundle. 

The most important aspect of this theory is that the moduli space of multi-monopoles $\M_Y(\mathfrak{s},\mathfrak{p})$ can be \emph{non-compact} for some $\mathfrak{p}$, unlike in the Seiberg--Witten case. This was explored by Haydys and Walpuski \cite{HaydysWalpuski2020}, who proved that a sequence of multi-monopoles either converges to another multi-monopole, or converges (after rescaling) to a $\Z/2$-harmonic spinor. These are two-valued harmonic spinors branching along a 1-dimensional subset of $Y$. Let $\mathcal{P}_Y$ be the space of parameters $\mathfrak{p}$, and let $\W_Y\subset\mathcal{P}_Y$ contain those parameters for which $\M_Y(\mathfrak{s},\mathfrak{p})$ is non-compact. A central conjecture in this subject is the following.

\begin{conj}\label{codim1}
    $\W_Y$ is a codimension one subset of $\mathcal{P}_Y$.
\end{conj}

The (signed) count of multi-monopoles $\#\M_Y(\mathfrak{s},\mathfrak{p})$ is constant in each connected component of $\mathcal{P}_Y\setminus\W_Y$, which is called a \emph{chamber}. This number can jump as $\mathfrak{p}$ crosses the \emph{wall} $\W_Y$, due to multi-monopoles degenerating to (or being born from) a $\Z/2$-harmonic spinor. This wall-crossing phenomenon means there is no hope of defining a topological invariant of $Y$ by counting multi-monopoles. Instead, the applications of these equations are of a rather different flavour. For example, wall-crossing can be used to prove the existence of $\Z/2$-harmonic spinors, which are difficult to construct by other means \cite{Doan-Walpuski2021}. Another more conjectural application is to define invariants of special holonomy manifolds using multi-monopoles \cite{Donaldson-Segal2009,DoanWalpuski_2019,Haydys2012}. In any applications of these equations, it will be crucial to understand examples of the wall-crossing phenomenon.

\subsection{Multi-vortices and wobbly bundles}\label{intro2}

The purpose of this article is to better understand the chamber structure and wall-crossing for $Y=S^1\times C$, the product of a circle and a closed Riemann surface $C$. Equipping $C$ with a volume form determines a metric $h_C$ and the product metric $h=dt^2+h_C$ on $Y$. One can also consider $\SU(2)$-connections $B$ and spin$^c$ structures pulled back from $C$. Such product-type parameters were studied by Doan, who showed that the multi-monopole equations on $Y$ reduce to the multi-vortex equations on $C$ \cite{Doan_2018}. Spin$^c$ structures on $C$ are indexed by an integer $d$, and an $\SU(2)$-connection on $C$ is equivalent to specifying a holomorphic $\SL(2,\C)$-bundle $E$ (rank $2$ with $\det=\O$). In terms of this complex-geometric data, Doan identified the moduli space of multi-vortices with
\begin{equation}\label{multivortex}
    \M_C^d(E)=\left\{(L,\a,\b)\Big|
    \begin{array}{l}
    \deg(L)=d,\ 0\neq\a\in H^0(L\otimes E),\\[2pt]
    \b\in H^0(L^*\otimes E^*\otimes K),\ \b\a=0
    \end{array}
    \right\}\Big/\sim.
\end{equation}

Here, $L$ is a degree $d$ holomorphic line bundle, $K$ is the canonical bundle, and the Serre duality pairing $\b\a$ takes values in $H^0(K)$. The equivalence relation is the action of automorphisms of $L$ on $\a$ and $\b$. Although for generic $E$, $\M_C^d(E)$ is a smooth projective variety of dimension $2d+1-g$, Doan showed that it is non-compact for $E$ in a subset of real codimension at least two \cite{Doan_2018}. This is higher than what Conjecture \ref{codim1} suggests, so wall-crossing does not occur in this setting and product-type parameters belong to just one chamber. Thakar has computed the count of multi-monopoles on $S^1\times C$ in this chamber \cite{thakar2025modulispacemultimonopolesriemann}:
\begin{equation}\label{productcount}
    \#\M_{S^1\times C}(\mathfrak{s}_d,\mathfrak{p}_{product})=\begin{cases}
        (-1)^{g-1}2^g{2g-2\choose 2d+1-g} & \textnormal{if }(g-1)/2\leq d\leq 3(g-1)/2,\\
        0 &\textnormal{otherwise.}
    \end{cases}
\end{equation}
Our first result describes this codimension two subset in terms of \emph{wobbly bundles}. A bundle $E$ is called wobbly if there is a non-zero nilpotent Higgs field $\phi:E\ra E\otimes K$. These were first studied by Drinfeld and Laumon in the context of the geometric Langlands correspondence, which relates certain local systems on $C$ to local systems on the complement of the wobbly locus $\W\subset\MM_C^{ss}$. Here, $\MM_C^{ss}$ is the the moduli space of holomorphic $\SL(2,\C)$-bundles. Pal and Pauly defined subvarieties $\W_k\subset\MM_C^{ss}$ for $1\leq k \leq g$ (see Section \ref{wobbly}) and proved that $\W=\W_{\lceil g/2\rceil}\cup\cdots\cup\W_g$ is a decomposition of $\W$ into irreducible divisors \cite{palpauly}.

\begin{theorem}\label{wobblystructure}
    Suppoes $C$ is a smooth curve of genus $g\geq 2$ and let $E$ be a semi-stable $\SL(2,\C)$-bundle on $C$. Then we have the following:
    \begin{enumerate}
        \item If $d<0$, then $\M_C^d(E)$ is empty. 
        \item If $0\leq d \leq 2g - 2$, then $\M_C^d(E)$ is non-compact if and only if
        $$E\in\W(g,d):=\bigcup_{k=1}^{g - |d+1-g|}\W_k.$$
        \item If $d>2g-2$, then $\M_C^d(E)$ is compact and biholomorphic to the projectivisation of the rank $2(d+1-g)$ vector bundle $V_E\ra J_C^d$ with fibre $H^0(L\otimes E)$ over $L$.
    \end{enumerate}
\end{theorem}

The proof uses that $\M_C^d(E)$ is non-compact if and only if it contains a point $(L,\alpha,\beta)$ with $\beta\neq 0$. Then the key idea is that $\beta\alpha=0$ implies $\phi=\alpha\beta$ is a non-zero nilpotent Higgs field. 

\subsection{Counting multi-monopoles on $S^1\times C$ using monodromy}\label{intro3}

To explore other chambers in $\mathcal{P}_{S^1\times C}$, we will continue using the product metric, but allow for non-constant paths of holomorphic $\SL(2,\C)$-bundles $E_t$ (equivalently, $\SU(2)$-connections $B_t$) satisfying $E_t=E_{t+1}$. Such a loop determines an $\SU(2)$-connection $B$ on $S^1\times C$. The K\"{u}nneth decomposition identifies every spin$^c$ structure on $S^1\times C$ with a pair $(d,\gamma)\in\Z\oplus H^1(C;\Z)$, with $(0,0)$ being the pullback of the canonical spin$^c$ structure on $C$. The integer $d$ is the degree of the spinor bundle on $\{t\}\times C$. This is a special case of the setting considered in \cite{wilson}, which deals more generally with mapping tori and non-product metrics. In what follows, we will need to consider only \emph{admissible} loops.

\begin{defn}\label{admissibleloops}
    For a given $d\in\Z$, a loop of holomorphic $\SL(2,\C)$-bundles $E_t$ is called admissible if $\M_C^d(E_t)$ is compact for all $t$.
\end{defn}

Smoothness of $\M_C^d(E)$ is equivalent to compactness \cite{Doan_2018}, so admissible loops determine a smooth fibre bundle of moduli spaces $\mathbf{M}\ra S^1$. It is equipped with an Ehresmann connection and there is an induced monodromy map $\Upsilon$ \cite{Salamon1999,wilson}. Define the rescaled metric $g_\e=dt^2+\e^2g_C$ and let $\mathfrak{p}_\e=(g_\e,B)\in\mathcal{P}_{S^1\times C}$. Then in the \emph{adiabatic limit} $\e\ra 0$, the multi-monopole equations for $\mathfrak{p}_\e$ formally degenerate to the equation for a parallel section of $\mathbf{M}$. This idea goes back to Salamon \cite{Salamon1999}, who described the adiabatic limit of the classical Seiberg--Witten equations in this way. The expectation is that for $\e>0$ sufficiently small, counting multi-monopoles is equivalent to counting parallel sections of $\mathbf{M}$, which in turn is equivalent to counting fixed points of $\Upsilon$. 

Before stating the main result of \cite{wilson}, which proves one half of this adiabatic correspondence, we define the subset $\Fix_{\mathfrak{s}}(\Upsilon)\subset\Fix(\Upsilon)$. Given $x\in\Fix(\Upsilon)$, there is a corresponding parallel section $s_x$ of $\mathbf{M}\ra S^1$. Composing with the forgetful map $\pi:\M_C^d(E)\ra J_C^d$ to the Jacobian torus, $s_x$ descends to $\pi(s_x):S^1\ra J_C^d$. This determines a class $\gamma_x\in H^1(C;\Z)\simeq \pi_1(J_C^d)$, so for $\mathfrak{s}=(d,\gamma)$ we can define 
$$\Fix_{\mathfrak{s}}(\Upsilon) = \{x\in\Fix(\Upsilon) : \gamma_x=\gamma\}.$$

\begin{theorem}[\cite{wilson}]\label{wilson}
    Let $E_t$ be a generic admissible loop. For $\e>0$, let $\mathfrak{p}_{\e}=(g_\e,B)$ be as defined above. Then for every spin$^c$ structure $\mathfrak{s}=(d,\gamma)$ and all sufficiently small $\e>0$, the fixed points of $\Upsilon$ are all non-degenerate and there is an injective map
    \begin{equation}\label{injection}
        \Fix_{\mathfrak{s}}(\Upsilon)\hookrightarrow \M_{S^1\times C}(\mathfrak{s},\mathfrak{p}_\e).
    \end{equation}
\end{theorem}

\begin{conj}\label{bijection}
    The map \eqref{injection} is a bijection and the orientation of each multi-monopole agrees with the index of the corresponding fixed point. That is,
    \begin{equation}\label{equality}
        \#\Fix_{\mathfrak{s}}(\Upsilon)=\#\M_{S^1\times C}(\mathfrak{s},\mathfrak{p}_\e).
    \end{equation}
\end{conj}

A homotopy of the loop $E_t$ varies the parameter $\mathfrak{p}_\e\in\mathcal{P}_{S^1\times C}$, and the fixed point count \eqref{equality} does not change as long as the compactness assumption holds along the homotopy. The work of Atiyah--Bott shows that the unstable strata are high codimension in the space of $\SU(2)$-connections \cite{atiyahbott}, so without loss of generality we can assume $E_t$ is semi-stable (see Section \ref{absection}). Theorem \ref{wobblystructure} then shows that the compactness assumption is fulfilled by loops $[E_t]\in\MM_C^{ss}$ avoiding the divisor $\W(g,d)$. Wall-crossing on $S^1\times C$ is expected to occur along homotopies passing through $\W(g,d)$, which can change the topology of $\mathbf{M}$ and the fixed point count of $\Upsilon$. By studying $\Upsilon$ for various loops around $\W(g,d)$, we can produce many examples of this phenomenon, as demonstrated by the theorem below.

\begin{theorem}\label{wallcrossingexamples}
    Suppose $g\geq 2$. Then the following statements hold:
    \begin{enumerate}
        \item If $d\leq (g-1)/2$ or $d\geq 3(g-1)/2$, then for all spin$^c$ structures $\mathfrak{s}=(d,\gamma)$ and all admissible loops of bundles we have $\#\Fix_{\mathfrak{s}}(\Upsilon)=0$.
        \item If $g=2d+1$ is odd, $\M_C^d(E)$ can be identified with the maximal line sub-bundles of $E$ and consists of $2^{g}$ points for $E\in\M_C^{ss}\setminus\W_{d+1}$. The resulting $2^g$-sheeted cover of $\M_C^{ss}\setminus\W_{d+1}$ has full monodromy group $S_{2^g}$. In particular, for every integer $0\leq k\leq 2^g$, there is an admissible loop around with $\#\Fix(\Upsilon)=k$.
        \item Let $d=g-1$. Then for each of the $2^{2g}$ irreducible components of $\W_g$, there is an admissible loop passing around it with
        \begin{equation}\label{lefschetzdehntwist}
            \#\Fix(\Upsilon)=(-1)^{g-1}2^g{2g-2\choose g-1} - 2\delta,
        \end{equation}
        where $\delta\in\{0,1\}$ is the parity of $g$.
    \end{enumerate}
\end{theorem}

Of course, $\#\Fix(\Upsilon)$ is just the Lefschetz number of $\Upsilon$, which is the sum of $\#\Fix_{(d,\gamma)}(\Upsilon)$ over all $\gamma\in H^1(C;\Z)$. The first part of Theorem \ref{wallcrossingexamples} follows from Theorem \eqref{wobblystructure}, the fact that $\W_k$ has codimension $>1$ for $k<\lceil g/2\rceil$, and equation \eqref{productcount}. At least when $g$ is odd, parts 2 and 3 indicate that wall-crossing on $S^1\times C$ is incredibly rich, even though we have not stratified the count of fixed points by spin$^c$ structure. Their proofs require a more explicit understanding of the monodromy map. The general idea is to consider families of holomorphic $\SL(2,\C)$-bundles $\E\ra C\times\DD$ over the unit disk $\DD\subset\C$ such that:
\begin{enumerate}
    \item The holomorphic disk $u:\DD\ra\MM_C^{ss}$ given by $u(t)=[E_t]$ passes transversely through a generic point $[E_0]$ of the divisor $\W(g,d)$,
    \item For all $t\neq 0$, $u(t)\notin\W(g,d)$. In particular, $\E|_{\del\DD}$ is an admissible loop.
\end{enumerate}
Over $t\neq 0$ there is a family of multi-vortex moduli spaces, but non-compactness poses issues for extending over $t=0$. Instead, we consider the \emph{framed} multi-vortex moduli space
\begin{equation}\label{framedmultivortex}
    \N_C^d(E)=\{(L,\alpha) : \deg(L)=d, 0\neq\alpha\in H^0(L\otimes E)\}/\sim.
\end{equation}
Doan proved that $\N_C^d(E)$ is always compact, and whenever $\M_C^d(E)$ is compact too, they are both smooth and biholomorphic \cite{Doan_2018}. By Theorem \ref{wobblystructure}, we expect $\N_C^d(E)$ to be singular (or possibly obstructed) for $E\in\W(g,d)$. This means that for $\E$ as above, the family of framed multi-vortex moduli spaces $\N_C^d(\E)\ra\DD$ will have a singular fibre over $t=0$. It restricts over the boundary $S^1=\del\DD$ to the family $\mathbf{M}\ra S^1$, so $\Upsilon$ can be studied by understanding the singularities of $\N_C^d(E_0)$ and applying Picard--Lefschetz theory. In Section \ref{deformations}, we study the deformation theory of $\N_C^d$ and prove the following result about its singularities.

\begin{theorem}\label{isolatedsingularities}
    Suppose $\max\{2,\lceil g/2\rceil\}\leq k\leq g$. Then for a generic bundle $E\in\W_k$, $\N_C^d(E)$ is of the expected dimension $2d+1-g$ and has $2(g-k)\choose d+1-k$ isolated hypersurface singularities.
\end{theorem}

When $d=g-1$ and $k=g$, we prove that the unique singularity of $\N_C^{d}(E)$ (for generic $E\in\W_g$) is an ordinary double point by studying its second order deformation theory. This shows that the monodromy map is a Dehn twist, and part 3 of Theorem \ref{wallcrossingexamples} follows. Part 2 of Theorem \ref{wallcrossingexamples} relies on identifying points of $\N_C^d(E)$ with maximal line sub-bundles of $E$. We show that these spaces form a $2^g$-sheeted cover $T\ra \M_C^{ss}\setminus\W_{d+1}$. By studying $\N_C^d(E)$ for a generic point $E\in\W_{d+1}$, we find that there is a transposition in the monodromy group. When combined with a theorem of Ballico \cite{maximalsubbundles}, this proves that the monodromy group of $T$ is the full symmetric group $S_{2^g}$.

We conclude the paper with a detailed discussion of the genus 2 case. Narasimhan and Ramanan proved that $\MM_C^{ss}\simeq \P^3$ \cite{NRgenus2}, inside which $\W_1$ is a singular Kummer surface and $\W_2$ is the union of 16 hyperplanes, classically known as tropes. These are arranged in a $(6,16)$-configuration, which we describe from the perspective of $\MM_C^{ss}$ in Section \ref{616}. We use this to prove the following result.

\begin{theorem}\label{genus2theorem}
    Let $C$ be a smooth genus $2$ curve.
    \begin{enumerate}
        \item Suppose $d\neq1$. Then for all spin$^c$ structures $\mathfrak{s}=(d,\gamma)$ and all admissible loops of bundles, we have $\#Fix_{\mathfrak{s}}(\Upsilon)=0$.
        \item Suppose $d=1$. Then $\M_C^1(E)\simeq\N_C^1(E)$ is a smooth genus $5$ curve whenever $E\notin\W(2,1)$. Moreover:
        \begin{enumerate}
            \item For an admissible loop passing once around the Kummer surface $\W_1$, the monodromy map is a composition of two commuting Dehn twists.
            \item For an admissible loop passing once around a trope hyperplane in $\W_2$, the monodromy map is a single Dehn twist.
        \end{enumerate}
    \end{enumerate}
\end{theorem}

Parts 2a and 2b of Theorem \ref{genus2theorem} show that there is a wealth of monodromy for $d=1$, but in terms of fixed point counts, we are unable to conclude anything more than $\#\Fix=-8$ for these loops. However, this does not preclude the existence of wall-crossing. For example, by understanding more about the vanishing cycles of these Dehn twists, it may be possible to change the total count of fixed points by composing these loops. It may also be possible to directly understand the refined fixed point counts $\#\Fix_{\mathfrak{s}}$. 

\subsection{Acknowledgements}

I would like to thank my advisor Aleksander Doan for his encouragement with this project. The idea of wobbly bundles being related to multi-vortices was first explained to me by Ollie Thakar and Daniel Santiago-Alvarez during a visit to Harvard in March 2026. This paper would not have been possible without their suggestion. Finally, many thanks to Bogdan Simeonov, Calum Crossley, and Yifan Zhao for patiently answering my questions about algebraic geometry. This work was supported by the Engineering and Physical Sciences Research Council [EP/S021590/1]. The EPSRC Centre for Doctoral Training in Geometry and Number Theory (The London School of Geometry and Number Theory), University College London.

\subsection{AI Disclosure}

Claude Opus 5.5 was used to find relevant papers in the literature and to assist with the details of some proofs. All ideas are the author's own and LLMs were not used to write any part of the final document. 

\section{Wobbly bundles and non-compactness}\label{wobbly}
We begin this section by reviewing the geometry of the Hitchin fibration and wobbly bundles. The main result relates wobbly bundles to the existence of Fueter sections, which are the source of non-compactness phenomena in multi-spinor Seiberg--Witten theory. Using Pal and Pauly's stratification of the wobbly locus \cite{palpauly}, we prove that the multi-vortex moduli space is non-compact whenever the auxilliary bundle lies in a certain collection of the strata. 

\subsection{The Hitchin fibration and wobbly bundles}

Let $C$ be a smooth projective curve of genus $g\geq 2$, with $K$ its canonical bundle. Denote by $\M_C^s$ (resp. $\MM_C^{ss}$) the moduli space of stable (resp. semi-stable) rank 2 vector bundles on $C$ with trivial determinant. The cotangent space $T^*\M_C^s$ can be partially compactified by the moduli space $\mathcal{H}_C^{ss}$ of semi-stable $\SL(2,\C)$-Higgs bundles. The Hitchin map is defined as
$$h:\mathcal{H}_C^{ss}\ra H^0(K^2), \hspace{3mm}(E,\phi)\mapsto \det\phi.$$
Hitchin proved that $\mathcal{H}_C^{ss}$ is an algebraic symplectic manifold (extending the cotangent symplectic structure on $T^*\M_C^s$), $h$ is proper, and its fibres are Lagrangian subvarieties of $\mathcal{H}_C^{ss}$ \cite{hitchin}. In particular, the nilpotent cone $h^{-1}(0)$ is Lagrangian. It decomposes as
$$h^{-1}(0)=\MM_C^{ss}\cup N,$$
where $\MM_C^{ss}$ corresponds to the Higgs bundles of the form $(E,0)$, and $N$ consists of Higgs bundles with a non-zero nilpotent Higgs field. Define the wobbly locus $\W\subset\MM_C^{ss}$ to be the image of $N$ under the rational forgetful map
$$p:\mathcal{H}_C^{ss}\dashrightarrow \MM_C^{ss},\hspace{3mm}(E,\phi)\mapsto E.$$
This is a closed subset consisting of wobbly bundles; those semi-stable bundles admitting a non-zero nilpotent Higgs field. The complement of the wobbly locus, which consists of \emph{very stable} bundles, was originally studied by Drinfeld and Laumon as an aspect of the geometric Langlands correspondence \cite{Laumon1988}, which has recently been proven \cite{GLC}. For $\SL(2,\C)$, this implies that certain $\textnormal{PSL}(2,\C)$-local systems on $C$ correspond to local systems on $\MM_C^{ss}\setminus\W$.

Drinfeld conjectured that $\W$ is a divisor. In the rank 2 case, this was proven by Pal and Pauly, who also described its irreducible components \cite{palpauly}. To state their result, define $\W_k^0$ to be the locus of all semi-stable bundles arising as extensions
\begin{equation}\label{wobblyseq}
    0\ra M\ra E\ra M^{-1}\ra 0,
\end{equation}
with $\deg(M)=1-k$ and $h^0(M^2\otimes K)>0$. Define $\W_k$ to be the closure of $\W_k^0$ in $\MM_C^{ss}$. 

\begin{remark}\label{wobblyequivalence}
    Bundles $E\in\W_k$ are wobbly, since for any non-zero map $M^{-1}\ra M\otimes K$ the composition $E\ra M^{-1}\ra M\otimes K\ra E\otimes K$ is a non-zero nilpotent Higgs field. The converse is also true, which can be shown by taking $M=\ker(\phi)\subset E$ (see \cite[Lemma 2.4]{dinh} for a full proof).
\end{remark}

\begin{theorem}\cite{palpauly}\label{ppthm}
    The wobbly locus $\W\subset \MM_C^{ss}$ decomposes as
    $$\W=\W_{\lceil g/2\rceil}\cup\dots\cup\W_g.$$
    Moreover, all of the $\W_k$ are irreducible divisors, except for $\W_g$ which has $2^{2g}$ irreducible components. 
\end{theorem}

\begin{remark}\label{lowerwobblycomponents}
    Pal and Pauly also proved that $\W_k$ is contained in $\W_{\lceil g/2\rceil}$ for $k<\lceil g/2\rceil$. Also, $\W_1$ is precisely the strictly semi-stable locus, which is also the singular locus of $\MM_C^{ss}$ when $g>2$.
\end{remark}

\begin{remark}\label{wobblyplanes}
    The $2^{2g}$ irreducible components of $\W_g$ correspond to the $2^{2g}$ line bundles $L\in J_C^{g-1}$ satisfying $L^2\simeq K$. Explicitly, the components are
    $$H_L = \{E\in\MM_C^{ss} : h^0(L\otimes E)>0\}.$$
    These divisors $H_L$ are all linearly equivalent and the associated line bundle generates the Picard group of $\MM_C^{ss}$.
\end{remark}

Another important result about the geometry of $\W$ was recently proven by Nollau.

\begin{theorem}\cite{nollau}\label{nollauthm}
    For $2\leq k \leq g$, suppose $E\in\W_k\cap\M_C^{s}$ is chosen generically. Then there is a unique line sub-bundle $M\subset E$ with $h^0(M^2\otimes K)>0$. Moreover, this $M$ has $h^0(M^2\otimes K)=1$. 
\end{theorem}

Note that $k=1$ is excluded because $\W_1$ is the strictly semi-stable locus. The upshot of having $h^0(M^2\otimes K)=1$ is that for generic $E$, there is a unique non-zero nilpotent Higgs field up to scaling. The following lemma relates this to the geometry of $\W$.

\begin{lemma}\label{wobblytangent}
    For $\max\{2,\lceil g/2\rceil\}\leq k\leq g$, suppose $E\in\W_k\cap \M_C^s$ is both a smooth point of $\W_k$ and generic in the sense of Theorem \ref{nollauthm}. Then the unique non-zero nilpotent Higgs field $\phi\in H^0(\End_0(E)\otimes K)$ satisfies
    $$T_E\W = \ker(\phi^*)\subset H^1(\End_0(E)) = T_E\MM_C^{ss}.$$
    Here, $\phi^*:H^1(\End_0(E))\ra \C$ is the Serre dual of $\phi$. 
\end{lemma}

\begin{proof}
    Consider the intersection 
    $$N':=N\cap p^{-1}(\M_C^s)=\{(E,\phi)\in\mathcal{H}_C^{ss} : E \textnormal{ stable, } \phi\neq 0\textnormal{ nilpotent}\},$$
    which is contained in $T^*\M_C^s\subset \mathcal{H}_C^{ss}$. By Theorem \ref{nollauthm}, the restricted projection map
    $$p:N'\ra \W\cap \M_C^s$$
    has 1-dimensional fibre $\C\phi$ over a generic smooth point $E\in\W_k\cap\M_C^s$, where $\phi$ is the unique non-zero nilpotent Higgs field associated to $E$. Suppose $v\in T_E(\W_k\cap \M_C^s)\subset H^1(\End_0(E))$, which determines a tangent vector $(v,0)$ in
    $$T_{(E,\phi)}N'\subset T_{(E,\phi)}T^*\M_C^s = H^1(\End_0(E))\oplus H^1(\End_0(E))^*.$$
    On the other hand, the tangent vector running along the fibre of $N'$ over $E$ is $(0,\phi^*)\in T_{(E,\phi)}N'$. Since a neighbourhood of $(E,\phi)$ in $N$ (and hence $N'$) is a smooth Lagrangian submanifold of $T^*\M_C^s$, we have
    $$\phi^*(v)=\omega((0,\phi^*),(v,0))=0.$$
    Then $T_E(\W_k\cap\M_C^s)\subset\ker(\phi^*)$, but for $\lceil g/2\rceil\leq k\leq g$ the inclusion is an equality as both of these spaces are of codimension 1 in $H^1(\End_0(E))$.
\end{proof}

\subsection{Multi-vortices and the existence of Fueter sections}

Let $d\in\Z$ and consider the multi-vortex moduli spaces $\M_C^d(E)$ and $\N_C^d(E)$ defined in Section \ref{intro}, where $E$ is a holomorphic $\SL(2,\C)$-bundle. The following result describes the generic behaviour of these moduli spaces.

\begin{theorem}\cite{Doan_2018}\label{genericsmoothness}
    Let $C$ be a smooth projective curve and let $d\in\Z$. Then for a semi-stable rank 2 bundle $E$ on $C$ with trivial determinant, the following conditions are equivalent:
    \begin{itemize}
        \item $\N_C^d(E)$ is a smooth projective variety. 
        \item $\M_C^d(E)$ is compact and biholomorphic to $\N_C^d(E)$.
        \item There is no triple $(L,\alpha,\beta)$ consisting of $L\in J_C^d$ and non-zero holomorphic sections $\alpha\in H^0(L\otimes E)$, $\beta\in H^0(L^*\otimes E^*\otimes K)$ satisfying $\beta\alpha=0\in H^0(K)$.
    \end{itemize}
    Moreover, these conditions hold for generic such bundles $E$. 
\end{theorem}

\begin{remark}\label{smoothnessremark}
    Before Doan's differential geometric proof of Theorem \ref{genericsmoothness}, similar results were already known to algebraic geometers. In particular, it was known that $\N_C^d(E)$ is the \emph{canonical blowup} \cite[Chapter II]{acgh} of the twisted Brill--Noether locus
    $$B_{1,d}^1(E)=\{L\in J^d_C : h^0(L\otimes E)>0\},$$
    giving it the structure of a projective scheme. Its smoothness at $(L,\alpha)$ was shown to be equivalent to the injectivity of the Petri map \cite{hhn}
    \begin{equation}\label{petri}
        \mu:H^0(L^{*}\otimes E^*\otimes K)\ra H^0(K),\hspace{3mm}\beta\mapsto \beta\alpha.
    \end{equation}
    Moreover, Teixidor i Bigas proved that this map is injective for generic $E$ and $C$ \cite{tib}. Theorem \ref{genericsmoothness} shows that the genericity assumption on $C$ is unnecessary. 
\end{remark}

In multi-spinor Seiberg--Witten theory on a 3-manifold, non-compactness of the moduli space is caused by the existence of Fueter sections \cite{HaydysWalpuski2020}. In the dimensional reduction to $C$, Doan gave the following complex-geometric description of these objects \cite{Doan_2018}. 

\begin{defn}
    Let $E$ be a rank 2 vector bundle on $C$ with trivial determinant. A Fueter section is a tuple $(L,\tilde{\alpha},\tilde{\beta},D_1,D_2)$ consisting of $L\in J_C^{d}$, divisors $D_1,D_2\geq 0$ such that $D_1+D_2$ is effective, and non-zero holomorphic maps $\tilde{\alpha},\tilde{\beta}$ forming the short exact sequence
    \begin{equation}\label{fueterseq}
        0\ra L^{-1}(D_1)\xrightarrow{\tilde{\alpha}} E\xrightarrow{\tilde{\beta}}L^{-1}K(-D_2)\ra 0.
    \end{equation}
\end{defn}

These objects correspond to the usual notion of a Fueter section on the 3-manifold $S^1\times C$, with singular set $S^1\times (D_1\cup D_2)$. Not every Fueter section appears as the limit of a sequence of multi-monopoles, but those that do are called limiting configurations \cite{Doan_2018}. 

\begin{defn}
    A limiting configuration is a Fueter section $(L,\tilde{\alpha},\tilde{\beta},D_1,D_2)$ where the divisors $D_1$ and $D_2$ are effective. 
\end{defn}

To prove Theorem \ref{wobblystructure}, we relate the exact sequences \eqref{wobblyseq} and \eqref{fueterseq}.

\begin{proof}[Proof of Theorem \ref{wobblystructure}]
    Firstly, suppose $d<0$. Then if there is a point $(L,\alpha,\beta)\in\M_C^d(E)$, the non-zero holomorphic map $\alpha:L^{-1}\ra E$ contradicts the semi-stabiliy of $E$. 

    Now fix $0\leq d\leq 2g-2$ and suppose $\M_C^d(E)$ is non-compact, so there is $(L,\a,\b)\in\M_C^d(E)$ with $\beta\neq 0$. Let $D_1$ and $D_2$ be the zero-divisors of $\a$ and $\b$ respectively. Then there are holomorphic maps $\tilde{\a}$ and $\tilde{\b}$ such that the following diagram commutes and is exact at $E$:
    \[\begin{tikzcd}
	    & {L^{-1}} &&& \\
	    0 & {L^{-1}(D_1)} & E & {L^{-1}K(-D_2)} & 0 \\
	    &&& {L^{-1}K}
	    \arrow[from=1-2, to=2-2]
	    \arrow["\alpha", from=1-2, to=2-3]
	    \arrow[from=2-1, to=2-2]
	    \arrow["{\tilde{\alpha}}"', from=2-2, to=2-3]
	    \arrow["{\tilde{\beta}}", from=2-3, to=2-4]
	    \arrow["\beta"', from=2-3, to=3-4]
	    \arrow[from=2-4, to=2-5]
	    \arrow[from=2-4, to=3-4]
    \end{tikzcd}\]
    Setting $M=L^{-1}(D_1)$, we find that $h^0(M^2\otimes K)=h^0(D_1+D_2)>0$ because of the condition $\det(E)=\O$. This means that $E\in\W_k$, where 
    $$1-k=\deg(M)=\deg(D_1)-d.$$
    Semi-stability of $E$ implies that $k\geq 1$. It follows from $\O(D_1+D_2)=M^{2}\otimes K$ that
    $$k\leq k + \deg(D_2)=k + 2(g-k)-\deg(D_1) = 2g-d-1.$$
    Combining this with $k\leq d+1$, we have
    $$k\leq \min\{d+1,2g-d-1\} = g-|d+1-g|.$$
    
    Conversely, suppose that $E\in\W_k$ for some $1\leq k \leq g-|d+1-g|$. Then there is line bundle $M$ of degree $1-k$ such that
    $$0\ra M\ra E\ra M^{-1}\ra 0$$
    is exact and there is a non-zero holomorphic section $\theta\in H^0(M^2\otimes K)$. Let $D$ be the zero-divisor of $\theta$, which has degree $2(g-k)\geq 2|d+1-g|$ and satisfies $\O(D)=M^2\otimes K$. Now choose any effective divisor $0\leq D_1<D$ of degree $d+1-k$, and let $D_2=D-D_1>0$ (which has degree $2g-k-d-1$). Define $L=M^{-1}(D_1)$, which has degree $d$. Moreover, we can define
    $$\alpha=\tilde{\alpha}\circ s_{1}:L^{-1}\ra E,\hspace{4mm}\beta=s_2^*\circ\tilde{\beta}:E\ra L^{-1}K,$$
    where $\tilde{\alpha},\tilde{\beta}$ are the maps appearing in the short exact sequence above, and $s_i\in H^0(\O(D_i))$ are non-zero sections. These are both non-zero and satisfy $\beta\alpha=0$, so $\M_C^d(E)$ is non-compact.

    As for the third point, $d>2g-2$ is large enough that $h^1(E\otimes L)=0$ for all $L\in J_C^d$. Then Riemann--Roch shows that the spaces $H^0(E\otimes L)$ assemble to a vector bundle $V_E$. We also find that $(L,\alpha,\beta)\in\M_C^d(E)$ necessarily has $\beta=0$, so it is compact. The map $\M_C^d(E)\ra J_C^d$ forgetting $\alpha,\beta$ is precisely the projection $\P V_E\ra J_C^d$, since $\alpha$ varies in $H^0(L\otimes E)$. 
\end{proof}

\begin{remark}
    When $d=2g-2$ and $E$ is stable, it is also the case that $\M_C^d(E)$ is the projectivisation of a vector bundle over $J_C^d$. This is consistent with the fact that $\W(g,2g-2)=\W_1$ is the strictly semi-stable locus. 
\end{remark}

\subsection{Unstable bundles and the Atiyah--Bott stratification}\label{absection}

So far, we have only considered semi-stable bundles. The following example shows that non-compactness can also occur when the bundle $E$ is unstable.

\begin{ex}
    Suppose $g\geq 1$ so that the canonical divisor $K$ is effective. Then $K+2x$ is effective for any $x\in C$, and we can choose two effective divisors $D_1,D_2\geq 0$ such that $K=D_1+D_2-2x\geq 0$. For $L=\O(D_1-x)$ and $E=\O(x)\oplus\O(-x)$ (which is unstable), the exact sequence
    $$0\ra \O(x)\xrightarrow{\tilde{\a}}\O(x)\oplus\O(-x)\xrightarrow{\tilde{\b}}\O(-x)\ra 0$$
    defines a limiting configuration since $D_1,D_2$ are effective.
\end{ex}  

While such unstable examples do exist, they are inconsequential to the monodromy picture described in Section \ref{intro}. To understand why, identify the (infinite-dimensional) space of holomorphic $\SL(2,\C)$-bundles with the space of $\SU(2)$-connections $\A$. Atiyah--Bott introduced a gauge-invariant stratification of $\A$, where the strata $\A_\mu$ are indexed by a ``slope vector'' $\mu$. In the case at hand, the slope vector is either $\mu=0$ or $\mu=(m,-m)$, where $m>0$ is an integer specifying the degree of the first term in the Harder--Narasimhan filtration of the holomorphic bundle $E_B$ induced by $B\in\A_\mu$. The stratum $\A_0$ consists of those connections inducing semi-stable holomorphic bundles, and all other strata correspond to unstable bundles. 

\begin{prop}[{\cite[p.~569]{atiyahbott}}]\label{atiyahbottcodim}
    The complex codimension of $\A_{(m,-m)}$ in $\A$ is equal to $2m+g-1$. 
\end{prop}

This comes from applying \cite[Equation 7.15]{atiyahbott} in our situation. It implies that for $g\geq 1$, the codimension of the unstable stratum determined by $m>0$ is 
\begin{equation}\label{unstablecodim}
    \codim_{\C}\A_{(m,-m)} = 2m+g-1\geq 2.
\end{equation}
Let $\widetilde{\W}(g,d)\subset\A_0$ be the union of the gauge orbits of connections $B\in\A$ with $[E_B]\in\W(g,d)$. By Proposition \ref{wobblystructure}, this consists of connections in $\A_0$ for which the multi-vortex moduli space is non-compact. The following corollary summarises our results. 

\begin{cor}\label{atiyahbottcorollary}
    Let $\A_{nc}\subset\A$ be the subset of $\SU(2)$-connections $B$ for which the multi-vortex moduli space $\M_C^d(E_B)$ is non-compact. Then
    $$\A_{nc}\subset \widetilde{\W}(g,d)\cup\bigcup_{_{m\geq 0}}\A_{(m,-m)},\hspace{4mm}\A_{nc}\cap\A_0=\widetilde{\W}(g,d),$$
    and any admissible loop in $\A\setminus\A_{nc}$ is homotopic to one contained in $\A_0\setminus\widetilde{\W}(g,d)$.
\end{cor}
This means that we can safely ignore the higher codimension unstable strata, and only ever need to compute the monodromy maps arising from admissible loops of semi-stable bundles.

\section{Deformation theory of multi-vortices}\label{deformations}
It follows from Theorem \ref{genericsmoothness} that $\M_C^d(E)$ is non-compact whenever $\N_C^d(E)$ is either singular (as a scheme) or of higher than expected dimension. In this section, we study the deformation theory of the moduli space of framed multi-vortices $\N_C^d(E)$ and families thereof. By describing the obstruction for a first order deformation to extend to second order, we can obtain some information about the singularities. 

\subsection{The deformation complex}

Let $\E\ra C\times S$ be a holomorphic family of rank 2 vector bundles on a smooth curve $C$ of genus $g\geq 2$. Then there is an associated family of (framed) multi-vortex moduli spaces $\N_C^d(\E)\ra S$. This is a special case of the more general space $G^k(\mathcal{V},\E)$ considered by Hitching, Hoff, and Newstead \cite[Section 3.4]{hhn}, which has the structure of a projective $S$-scheme. In our case, $k=1$ and $\mathcal{V}$ is a Poincar\`{e} line bundle on $C\times J^d_C$. They also studied conditions for a point $(L,\alpha,s)\in\N_C^d(\E)$ to be smooth. To state their result, define $\mu$ by
\begin{equation}\label{mu}
    \mu:H^0(L^{-1}\otimes E_s^*\otimes K)\ra T_s^*S\oplus H^0(K),\hspace{3mm}\beta\mapsto (\kappa^*(\phi),\beta\alpha),
\end{equation}
where $\phi=\alpha\beta\in H^0(\End(E_s)\otimes K)$ and $\kappa^*:H^0(\End(E_s)\otimes K)\ra T_s^*S$ is the adjoint (with respect to the Serre duality pairing) of the Kodaira--Spencer map $\kappa:T_sS\ra H^1(\End(E_s))$ of $\E$ at $s$. The adjoint map $\mu^*$ is
$$\mu^*:T_sS\oplus H^1(\O)\ra H^1(L\otimes E_s),\hspace{3mm}(u,v)\mapsto\kappa(u)\cup\alpha + v\cup\alpha.$$

\begin{prop}\cite{hhn}\label{familydef}
    Suppose $(L,\alpha,s)\in\N_C^d(\E)$. Then if $d\pi$ is the derivative of the projection $\pi:\N_C^d(\E)\ra S\times J^d_C$, for every $s\in S$ there is an exact sequence 
    \begin{equation}\label{familydefcomplex}
        0\ra H^0(L\otimes E_s)/\C\alpha\ra T_{(L,\alpha,s)}\N_C^d(\E)\xrightarrow{d\pi}T_sS\oplus H^1(\O)\xrightarrow{\mu^*}H^1(L\otimes E_s)
    \end{equation}
    for the Zariski tangent space. Moreover, $\N_C^d(\E)$ is smooth and of the expected dimension at $(L,\alpha,s)$ if and only if $\mu^*$ is surjective (or equivalently, $\mu$ is injective). 
\end{prop}

In our examples, the parameter space $S$ will be the unit disk $\DD\subset\C$. We will typically choose the central fibre of a family $\E\ra C\times\DD$ to be a wobbly bundle, so that the central fibre of $\N_C^d(\E)$ is singular. In what follows, we describe more carefully the families $\E$ we are interested in, and prove that the total space of $\N_C^d(\E)$ is smooth. Recall that any family $\E\ra C\times S$ determines a morphism $u:S\ra \MM_C^{ss}$ given by $u(s)=[\E|_{C\times\{s\}}]$.

\begin{defn}\label{transversefamily}
    A family of semi-stable holomorphic $\SL(2,\C)$-bundles $\E\ra C\times\DD$ is \emph{transverse to }$\W_k$ if the associated holomorphic disk $u:\DD\ra\MM_C^{ss}$ intersects $\W_k$ transversely at $u(0)$. Such a family will be called \emph{small} if $u(t)\notin\W(g,d)$ for all $t\neq 0$. 
\end{defn}

\begin{remark}
    The admissible loop $\E|_{\del\DD}$ obtained from a small family $\E\ra C\times\DD$ is homotopic to a loop confined to a neighbourhood of $E_0\in\W_k$. 
\end{remark}

\begin{lemma}\label{stablefamily}
    Suppose $u:\DD\ra\M_C^{s}$ is a holomorphic disk. Then up to isomorphism, there is a unique family of stable bundles $\E\ra C\times\DD$ such that $\E|_{C\times\{t\}}\simeq u(t)$. 
\end{lemma}

\begin{proof}
    There is a universal projective bundle $\P\ra C\times\M_C^s$ with $\P|_{(x,E)}$ canonically identified with $\P E_x$. Lifting $\P$ to a universal vector bundle on $\M_C^s$ is obstructed by its Brauer class $\beta$, which generates the Brauer group $H^2(\M_C^s,\O^\times)$ \cite{bbgn}. Since $\gcd(2,0)=2$ in our case, the Brauer group is $\Z/2\Z$ and $\beta$ is non-trivial. However, its pullback under $u:\DD\ra\M_C^s$ is trivial since $H^2(\DD,\O^\times)$ vanishes. This means that $\P$ lifts along the image of $u$ to define the required family, which is unique up to isomorphism because stable bundles are simple.
\end{proof}

\begin{prop}\label{stablefamilysmoothness}
    Suppose $\max\{2,\lceil g/2\rceil\}\leq k\leq g$ and let $E_0$ be a generic smooth point of $\W_k$. Then there exists a small family of stable bundles $\E\ra C\times\DD$ passing transversally through $\W_k$ at $E_0$. Moreover, the total space of $\N_C^d(\E)$ is smooth.
\end{prop}

\begin{proof}
    Since $\max\{2,\lceil g/2\rceil\}\leq k$, we can assume that $E_0$ is stable. In a local holomorphic chart around $E_0\in\M_C^s$, one can always find a holomorphic disk passing transversely through $\W_k$. Lemma \ref{stablefamily} then finds a family $\E$ passing transversely through $E_0$. Since $E_0$ is chosen generically, the holomorphic chart can be chosen small enough that it does not intersect any other $\W_{k'}$ for $k'\neq k$, so the family $\E$ is small. 

    Now let $(L,\alpha,t)\in \N_C^d(\E)$ and first suppose that $t\neq 0$. Since $E_t\notin\W(g,d)$, the fibre $\N_C^d(E_t)$ over $t$ is compact and smooth by Theorem \ref{wobblystructure}, so the only $\beta\in H^0(L^{-1}\otimes E_t^*\otimes K)$ with $\beta\alpha=0$ is $\beta=0$. This shows that the map $\mu$ is injective, so $\N_C^d(\E)$ is smooth for $t\neq 0$. Now suppose $t=0$, so there is $0\neq\beta\in H^0(L^{-1}\otimes E_0^*\otimes K)$ with $\beta\alpha=0$. Moreover, $\phi=\alpha\beta$ is a non-zero nilpotent Higgs field for $E_0$, which is unique up to scaling by our genericity assumption and Theorem \ref{nollauthm}. Lemma \ref{wobblytangent} implies that $\kappa(\del_t)\notin\ker\phi^*$, because otherwise, $\E$ deforms along the direction of $\W(g,d)$, contradicting transversality. This is equivalent to $\kappa^*(\phi)=dt\neq 0$, so $\mu$ is injective and $\N_C^d(\E)$ is also smooth for $t=0$. 
\end{proof}

\subsection{First order deformation theory}

For a fixed bundle $E$, more can be said about the deformation theory of $\N_C^d(E)$. A point $(L,\alpha)\in\N_C^d(E)$ determines an exact sequence of sheaves
\begin{equation}\label{alphasequence}
    0\ra L^{-1}\xrightarrow{\alpha}E\ra Q\ra 0,
\end{equation}
where $Q=E/\alpha(L^{-1})$. Specifying $(L,\alpha)$ is equivalent to specifying the quotient sheaf $Q$, so $\N_C^d(E)$ can be identified with the Quot scheme $\Quot^{1,d}(E)$ consisting of quotients $[E\twoheadrightarrow Q]$ of rank $1$ and degree $d$. Its Zariski tangent space at $[E\twoheadrightarrow Q]$ is $\Hom(L^{-1},Q)=H^0(L\otimes Q)$ \cite[Thm.~6.4.9]{fga}. In the lemma below, we spell out explicitly how this space parameterises first order deformations of $(L,\alpha)$. 

\begin{lemma}\label{firstorderdef}
    First order deformations of the pair $(L,\alpha)$ correspond to elements of $H^0(L\otimes Q)$. 
\end{lemma}

\begin{proof}
    A first order deformation of $(L,\alpha)$ is a line bundle $L_1\ra C\times\Spec\C[t]/(t^2)$ with $L_1|_C=L$, along with a section $\alpha_1$ of $L_1\otimes E$ with $\alpha_1\equiv\alpha\pmod{t^2}$. Let $\mathcal{U}=\{U_i\}_{i\in I}$ be a good cover of $C$ by affine open subsets, and let $U_{ij}=U_i\cap U_j$. The transition functions of $L_1\otimes E$ are of the form $g_{ij}(1+tf_{ij})$, where $g_{ij}\in \GL(L\otimes E)(U_{ij})$ are the transition functions of $L\otimes E$ and $f_{ij}\in \O(U_{ij})$. The section $\alpha_1$ is locally of the form $\alpha+ts_i$ over $U_i$ and must satisfy
    \begin{equation}\label{1storderdef}
        \alpha+ts_j=(1+tf_{ij})(\alpha+ts_i)
    \end{equation}
    to be globally well-defined. This amounts to $s_j-s_i=f_{ij}\alpha$, so if $q:L\otimes E\ra L\otimes Q$ is the quotient map, $q(s_j)=q(s_i)$ defines a global section of $L\otimes Q$. Conversely, $H^1(U_i,\O)=0$ since $U_i$ is affine, so any $x\in H^0(L\otimes Q)$ admits local lifts $s_i\in (L\otimes E)(U_i)$. Since $q(s_j-s_i)=0$, there is $f_{ij}\in\O(U_{ij})$ such that $s_j-s_i=f_{ij}\alpha$. This is equivalent to \eqref{1storderdef}, from which the first order deformation $(L_1,\alpha_1)$ can be recovered. 
\end{proof}

\begin{remark}\label{coboundaryremark}
    It is worth noting that $f_{ij}$ appearing in the proof of Lemma \ref{firstorderdef} defines a 1-cocycle in $H^1(\O)$. This is the evaluation of the coboundary map $\delta:H^0(L\otimes Q)\ra H^1(\O)$ at $x$. More explicitly, it is the cup product of $x$ with the extension class of $L\otimes E$. 
\end{remark}

It will also be useful to consider the long exact sequence associated with \eqref{alphasequence}:
\begin{equation}\label{deformationLES}
    0\ra H^0(L\otimes E)/\C\alpha\ra H^0(L\otimes Q)\xrightarrow{\delta}H^1(\O)\xrightarrow{\mu^*}H^1(L\otimes E)\ra H^1(L\otimes Q)\ra 0.
\end{equation}
The first four terms agree with \eqref{familydefcomplex} and Lemma \ref{firstorderdef} identifies $\delta$ with the derivative $d\pi$.

\subsection{Second order deformation theory}

To obtain more information about the singularities of $\N_C^d(E)$, we need to understand the second-order deformation theory around a point $(L,\alpha)$. The following lemma finds the obstruction to extending a first order deformation to a second order deformation. It has the same form as the obstruction appearing in \cite{nasu} for deformations of curves. 

\begin{lemma}\label{secondorderdef}
    Suppose a first order deformation of $(L,\alpha)\in\N_C^d(E)$ is determined by $x\in H^0(L\otimes Q)$, as described above. Then it extends to second order deformation if and only if
    $$\delta(x)\cup x\in H^1(L\otimes Q)=0,$$
    where $\delta:H^0(L\otimes Q)\ra H^1(\O)$ is the coboundary map induced by \eqref{alphasequence}.
\end{lemma}

\begin{proof}
    Given a first order deformation $x\in H^0(L\otimes Q)$, choose local lifts $s_i\in (L\otimes E)(U_i)$ satisfying \eqref{1storderdef} for some $f_{ij}$ on $U_{ij}$. Consider local extensions to a second order deformation $\alpha+ts_i+t^2r_i$ on $U_i$. These patch together to a global second order deformation if and only if there is $1+tf_{ij}+t^2h_{ij}$ such that
    \begin{equation}\label{2ndorderdef}
        \alpha+ts_j+t^2r_j=(1+tf_{ij}+t^2h_{ij})(\alpha+ts_i+t^2r_i).
    \end{equation}
    This amounts to the relation $r_j-r_i=h_{ij}\alpha + f_{ij}s_j$, which is equivalent to $q(r_j)-q(r_i)=f_{ij}x$. Since $q(r_j)-q(r_i)$ is a coboundary, it vanishes in $H^1(L\otimes Q)$. On the other hand, Remark \ref{coboundaryremark} shows that $f_{ij}x$ is a 1-cocycle representing $\delta(x)\cup x$, so this is the obstruction to extending $x$ to a second order deformation. 
\end{proof}

\begin{prop}\label{ordinarydoublepoint}
    Suppose $(L,\alpha)\in\N_C^d(E)$ has $h^1(L\otimes Q)=1$ and the quadratic form
    $$q:H^0(L\otimes Q)\ra H^1(L\otimes Q),\hspace{3mm}q(x)=\delta(x)\cup x$$
    is non-degenerate. Then $(L,\alpha)$ is an ordinary double point. 
\end{prop}

\begin{proof}
    Let $\N=\N_C^d(E)$ and $\mathfrak{m}\subset\O_{\N,z}$ the ideal of germs vanishing at $z=(L,\alpha)$. As shown in \cite[Thm.~6.4.9]{fga}, the local ring $\O_{\N,z}$ has an obstruction theory in the sense of \cite[Sect.~11]{hartshornedef}. The obstruction space is $\operatorname{Ext}^1(Q,L^{-1})=H^1(L\otimes Q)$ , which is one-dimensional by assumption. Let $n=\dim\mathfrak m/\mathfrak m^2=h^0(L\otimes Q)$ and choose local coordinate functions $y_1,\dots,y_n\in\mathfrak m$ whose classes form a basis of $\mathfrak m/\mathfrak m^2$. They define a surjection $\C\{y_1,\dots,y_n\}\to\O_{\N,z}$ with kernel $I$, corresponding to a closed embedding of germs $(\N,z)\hookrightarrow (\C^n,0)$. Here, $R:=\C\{y_1,\dots,y_n\}$ is the ring of convergent power series in the variables $y_i$, which has a unique maximal ideal $\mathfrak{m}_R=(y_1,\dots,y_n)$. That the $y_i$ form a basis shows that $\mathfrak m_R/\mathfrak m_R^2\to\mathfrak m/\mathfrak m^2$ is an isomorphism and $I\subset \mathfrak{m}_R^2$. Since $R$ is a regular local ring, \cite[Thm.~11.1.]{hartshornedef} shows that $I$ is generated by at most $h^1(L\otimes Q)=1$ elements. Then $I=(f)$ for some (possibly zero) $f\in\mathfrak m_R^2$ and $\O_{\N,z}\cong R/(f)$.

    Identify the Zariski tangent space $T_z\N$ with $\C^n$ using the basis dual to $y_1,\dots,y_n$. A tangent vector $x\in\C^n$ corresponds to the homomorphism $R/(f)\ra\C[t]/(t^2)$ given by $y_i\mapsto tx_i$. It extends to a second order deformation of $z$ if and only if this homomorphism lifts to $\C[t]/(t^3)$. Such a lift is given by $y_i\mapsto tx_i+t^2w_i$, where $w_i\in\C$ satisfy $f(tx+t^2w)\in (t^3)$. Using $f\in\mathfrak{m}_R^2$, write $f=f_2+f_3+\cdots$ as a sum of degree $k$ homogeneous polynomials $f_k$. Since $f_k(tx+t^2w)\in(t^3)$ for $k\geq 3$ and $f_2(tx+t^2w)=t^2f_2(x+tw)$, we have
    $$f(tx+t^2w)\equiv t^2f_2(x)\pmod{t^3},$$
    independently of $w$. This means that $x$ lifts to a second order deformation if and only if $f_2(x)=0$, so $\{f_2=0\}=\{q=0\}$ by Lemma \ref{secondorderdef}. Then there are $r>0$ and $g\in R$ such that $q=gf_2^r$. Both $f_2$ and $q$ are homogeneous of degree 2, so $g$ is constant, $r=1$, and $g=0$ is ruled out by non-degeneracy of $q$. Then $f_2$ is also non-degenerate and in particular, non-zero. Finally, the holomorphic Morse lemma finds local coordinates $u_i$ near $0\in\C^n$ in which $f=u_1^2+\cdots+u_n^2$, so $\O_{\N,z}\cong\C\{u_1,\dots,u_n\}/(u_1^2+\cdots+u_n^2)$ and $z$ is an ordinary double point. 
\end{proof}

\subsection{Isolated singularities}

The assumptions of Proposition \ref{ordinarydoublepoint} may be too much to hope for in general, particularly the non-degeneracy condition. We now prove Theorem \ref{isolatedsingularities}, which shows that the singularities are quite reasonable, at least generically. 

\begin{proof}[Proof of Theorem \ref{isolatedsingularities}]
    Suppose $E\in\W_k\cap \M_C^s$ is generic in the sense of Theorem \ref{nollauthm}, so that there is a unique bad line sub-bundle $M\subset E$ of degree $k-1$ satisfying $h^0(M^2\otimes K)=1$. Let $D$ be the (unique) effective divisor with $\O(D)=M^2\otimes K$, and consider the short exact sequence
    \begin{equation}\label{Msequence}
        0\ra M\xrightarrow{\tilde{\a}} E\xrightarrow{\tilde{\b}} M^{-1}\ra 0.
    \end{equation}
    The proof of Proposition \ref{wobblystructure} shows that every point $(L,\alpha,\beta)\in\M_C^d(E)$ with $\beta\neq 0$ solving $\beta\alpha=0$ arises in the following way. Choose an effective divisor $0\leq D_1\leq D$ of degree $d+1-k$, of which there are $2(g-k)\choose d+1-k$ choices. Also let $D_2=D-D_1\geq 0$. Then $L=M^{-1}(D_1)$ and $\a,\b$ must be of the form
    $$\a=\tilde{\a}\circ s_{D_1},\hspace{3mm}\b=s_{D_2}^*\circ\tilde{\b},$$
    where $s_{D_{i}}\in H^0(D_i)$. Since $0<h^0(D_i)\leq h^0(D)=1$, $\a$ and $\b$ are uniquely defined up to rescaling. There are then $2(g-k)\choose d+1-k$ isolated points at which $\N_C^d(E)$ is singular or of greater than expected dimension.

    We now prove that at any of these points, the obstruction space $H^1(L\otimes Q)$ is one-dimensional. Note that the injective map $\tilde{\a}:L^{-1}(D_1)\ra E$ descends to an injective map of sheaves $\O_{D_1}=L^{-1}(D_1)/L^{-1}\ra E/\alpha(L^{-1})=Q$. Its cokernel is $M^{-1}=KL^{-1}(-D_2)$ by the exactness of \eqref{Msequence}, so we obtain the exact sequence
    \begin{equation}\label{torsionsequence}
        0\ra \O_{D_1}\ra Q\ra KL^{-1}(-D_2)\ra 0.
    \end{equation}
    Tensoring by $L$ and taking cohomology, we get
    $$h^1(L\otimes Q)=h^1(K(-D_2))=h^0(D_2)=1,$$
    since $h^1(\O_{D_1})=0$. Now, the argument of Proposition \ref{ordinarydoublepoint} shows that $\N_C^d(E)$ is cut out by a single function $f$ near one of the isolated singular points. The possibility of $f=0$ is ruled out by nearby points being smooth and of the expected dimension. 
\end{proof}

\section{Monodromy around the wobbly divisor}
\label{monodromysection}

This section analyses two special cases where the assumptions of Proposition \ref{ordinarydoublepoint} can be verified, culminating in the proof of Theorem \ref{wallcrossingexamples}. 

\subsection{Monodromy around $\W_{d+1}$ for $g=2d+1$}

The study of maximal line sub-bundles of stable bundles on curves is a classical topic in algebraic geometry. Associated to a stable rank $2$ vector bundle on $C$ is the quantity
\begin{equation}\label{maximalquantity}
    s(E)=\deg(E)-2\max_{M\subset E}\{\deg(M)\},
\end{equation}
for which $0<s(E)\leq g$. An important result in the theory is that for $g=2d+1$ odd, every $E$ in an open dense subset of $\MM_C^{s}$ has $s(E)=g-1$ and there are $2^g$ line sub-bundles $M\subset E$ of the maximal degree $d=(g-1)/2$. This has many proofs, see for example \cite{lange85,Oxbury2000VarietiesOM,BIGAS_2004}. 

Notice that the space $M(E)$ of maximal line sub-bundles of $E$ can be identified with $\N_C^d(E)$. Indeed, when its expected dimension is zero (as is the case here), Thakar showed that $\N_C^{d}(E)$ consists of $2^g$ points for generic $E$ \cite{thakar2025modulispacemultimonopolesriemann}. Below, we use Theorem \ref{wobblystructure} to show that $\Omega$ can actually be taken as the complement of $\W(2d+1,d)=\W_{d+1}$. The idea that smoothness of $M(E)$ is related to the existence of non-zero nilpotent Higgs fields appears in \cite[Prop.~1.4]{Oxbury2000VarietiesOM}, although the later work of Pal and Pauly was required for an explicit understanding of the wobbly locus \cite{palpauly}. 

\begin{lemma}
    Let $g=2d+1$. Then there is a $2^g$-sheeted covering space $T\ra\Omega:=\M_C^{s}\setminus\W_{d+1}$ whose fibre over $E\in\Omega$ is $\N_C^d(E)\simeq M(E)$. Moreover, $\Omega$ is the maximal open subset of $\M_C^{s}$ over which such a covering space exists. 
\end{lemma}

\begin{proof}
    Let $\mathcal{V}\ra J_C^{-d}$ be the rank $2g-2$ vector bundle with fibre $\Ext^1(M^{-1},M)=H^1(M^{2})$ over $M\in J_C^{-d}$. There is a rational map
    $$\Phi:\P\mathcal{V}\dashrightarrow \M_C^{s}$$
    sending $(M,[e])$ to the bundle $E_e$ determined by the extension class $e$. It is well-defined on the dense open subset $\P\mathcal{V}^s$ consisting of pairs $(M,[e])$ for which $E_e$ is stable, and its fibres $\Phi^{-1}(E)$ can be identified with $M(E)$ and therefore $\N_C^d(E)$. Define
    $$\mathcal{Z}=\{M\in J_C^{-d} : h^0(M^2\otimes K)>0\}$$
    and note that $\Phi$ maps $\P\mathcal{V}^s|_{\mathcal{Z}}$ onto $\W_{d+1}\subset\M_C^{s}$. Define
    $$\Omega=\M_C^{s}\setminus\W_{d+1},\hspace{3mm}T=\P\mathcal{V}^s|_{J_C^{-d}\setminus\mathcal{Z}},$$
    which are both smooth. Theorems \ref{genericsmoothness} and \ref{wobblystructure} show that the fibres of $\Phi:T\ra\Omega$ are smooth and zero-dimensional. There is an open dense subset $\Omega'\subset\Omega$ such that $T|_{\Omega'}\ra\Omega'$ is a $2^{g}$-sheeted covering space, so $T\ra\Omega$ is too. That $\Omega$ is the maximal subset over which $\P\mathcal{V}^s$ restricts to a covering space follows from either \cite[Thm.~1.2]{palpauly} or Theorem \ref{wobblystructure}.
\end{proof}

We are interested in understanding the monodromy map $\pi_1(\Omega)\ra S_{2^g}$ of this covering space. Ballico proved the following result about its image $G$, the monodromy group. 

\begin{prop}[{\cite[Prop.~0.2]{maximalsubbundles}}]\label{twotransitivity}
    The monodromy group $G\subset S_{2^g}$ of the covering $T\ra\Omega$ is $2$-transitive. That is, for every two pairs $(x,y),(x',y')\in G\times G$, there exists $g\in G$ such that $x'=gx$ and $y'=gy$. 
\end{prop}

In the remainder of this section, we prove there is a transposition $\tau$ in $G$, obtained from a loop in $\Omega$ passing once around a generic point of $\W_{d+1}$. The 2-transitivity of $G$ then implies that $G=S_{2^g}$, since every other transposition $\tau'$ can be obtained from $\tau$ by conjugation. This method of proving that the monodromy group is the full symmetric group was originally applied by Harris to various enumerative problems in algebraic geometry \cite{harrisenumerative}.

\begin{lemma}\label{maximalsubbundledegeneration}
    Suppose $E\in\W_{d+1}$ is chosen generically, so that in particular, it has a unique bad line sub-bundle $\alpha:M\hookrightarrow E$ of degree $-d$. Then $\N_C^{d}(E)$ is a zero-dimensional scheme of length $2^g$. It contains $2^{g}-2$ reduced points and an isolated double point at $(M^{-1},\alpha)$.
\end{lemma}

\begin{proof}
    The unique bad line sub-bundle $M$ satisfies $h^0(M^2\otimes K)=1$, while any other $M'\not\simeq M$ has $h^0((M')^2\otimes K)=0$. It then follows from Serre duality and Riemann--Roch that $h^0(M^{-2})=h^1(M^{-2})=1$ and $h^0((M')^{-2})=h^1((M')^{-2})=0$. This shows that all points in $\N_C^d(E)$ are smooth except for $(M^{-1},\alpha)$ (see also \cite[Thm.~1.2]{palpauly}). This is consistent with Theorem \ref{isolatedsingularities} because we expect ${2(g-(d+1))\choose 0} = 1$ isolated singularities. To prove that $(M^{-1},\alpha)$ is a double point, by Proposition \ref{ordinarydoublepoint} it suffices to check that the quadratic form
    $$q:H^0(M^{-2})\ra H^1(M^{-2}),\hspace{3mm}x\mapsto \delta(x)\cup x$$
    is non-zero for generic $E\in\W_{d+1}$. For non-zero $\sigma\in H^0(M^2\otimes K)$ and $x\in H^0(M^{-2})$, we have
    $$\langle q(x),\sigma\rangle=\langle e, x\otimes x\otimes \sigma\rangle,$$
    where the brackets denote the Serre duality pairing and $e\in\Ext^1(M^{-1},M)=H^1(M^2)$ is the extension class of $E$. The pairing is non-degenerate and $x\otimes x\otimes\sigma\neq 0$, so $q(x)=0$ if and only if $e$ is contained in the hyperplane
    $$\Pi_M=(x\otimes x\otimes\sigma)^{\perp}\subset H^1(M^2).$$
    Define $\mathcal{Z}_1=\{M\in\mathcal{Z} : h^0(M^2\otimes K)=1\}$ and note that $\mathcal{Z}\setminus\mathcal{Z}_1$ has codimension at least 1 in $\mathcal{Z}$. The restriction $\P\mathcal{V}^s|_{\mathcal{Z}_1}$ is mapped into $\W_{d+1}$ by $\Phi$, and Theorem \ref{nollauthm} shows that its image is dense. Since
    $$\{(M,[e])\in\P\mathcal{V}^s:M\in\mathcal{Z}_1,\ [e]\in\Pi_M\}$$
    has codimension at least 1 inside $\P\mathcal{V}^s|_{\mathcal{Z}_1}$, its image under $\Phi$ is a subset of codimension at least 1 inside $\W_{d+1}$. So for a generic $E\in\W_{d+1}$, we can ensure both that there is a unique maximal line sub-bundle $M\subset E$ and that the quadratic map $q$ is non-zero at this $M$. Then $(M^{-1},\alpha)\in\N_C^d(E)$ is a double point, locally cut out by the equation $x^2=0$. It follows from \cite[Cor.~2.6]{gillam2011} that the zero-dimensional scheme $\N_C^d(E)$ has length $2^g$.
\end{proof}

Proposition \ref{stablefamilysmoothness} and Lemma \ref{maximalsubbundledegeneration} show that the restriction of $T$ to a holomorphic disk passing transversely through $\W_{d+1}$ is a branched covering of degree $2^g$ with one ramification point. This proves the following corollary. 

\begin{cor}\label{transposition}
    Let $\E\ra C\times\DD$ be a family of stable bundles passing transversely through a generic point $E_0\in\W_{d+1}$. Then the monodromy map of $\N_C^{d}(\E)|_{\del\DD}\ra \del\DD$ is a transposition swapping two of the $2^g$ points in a fibre and leaving the others fixed. 
\end{cor}

\begin{cor}\label{fullmonodromy}
    The monodromy group $G$ of the $2^g$-sheeted covering $T\ra \Omega$ is the full symmetric group $S_{2^g}$. In particular, this proves part $2$ of Theorem \ref{wallcrossingexamples}. 
\end{cor}

\subsection{Monodromy around $\W_g$ for $d=g-1$}

Recall from Remark \ref{wobblyplanes} that $\W_g$ consists of $2^{2g}$ irreducible components $H_L$, indexed by the line bundles $L$ such that $L^2=K$. 

\begin{prop}\label{Wg_odp}
    Suppose $L\in\Pic^{g-1}_C$ satisfies $L^2=K$ and choose a generic stable bundle $E\in H_L$. Then the moduli space $\N_C^{g-1}(E)$ has a unique ordinary double point at $(L,\alpha)$ and is otherwise smooth of the expected dimension.
\end{prop}

\begin{proof}
    By Theorem \ref{nollauthm}, a generic stable bundle $E\in H_L$ has $L^{-1}$ as its unique bad line sub-bundle, with $h^0(L\otimes E)=1$. This also precludes the possibility of $E\in H_{L'}$ for some $L'\neq L$ with $(L')^2=K$. Let $\alpha\in H^0(L\otimes E)$ be a non-zero section, which is unique up to scaling. Since ${2(g-k)\choose d+1-k} = {0\choose 0} = 1$, Proposition \ref{isolatedsingularities} shows that $\N_C^{g-1}(E)$ has a unique isolated singularity at $(L,\alpha)$. To see that it is an ordinary double point, note that $Q=E/\alpha(L^{-1})=L$ and the sequence \eqref{deformationLES} at $(L,\alpha)$ becomes
    $$0\ra H^0(K)\xrightarrow{\delta}H^1(\O)\xrightarrow{\mu^*}H^1(L\otimes E)\ra H^1(K)\ra 0.$$
    The obstruction space is $H^1(K)$, which is one-dimensional. Riemann--Roch gives us $h^1(M\otimes E)=1$, so the last map is an isomorphism and $\mu^*$ is zero. Then $\delta$ is an isomorphism, so the non-degeneracy of the Serre duality pairing $H^0(K)\otimes H^1(\O)\ra H^1(K)$ shows that the quadratic form $q(x)=\delta(x)\cup x$ is non-degenerate. Proposition \ref{ordinarydoublepoint} then shows that $(L,\alpha)$ is an ordinary double point.
\end{proof}

\begin{proof}[Proof of Theorem \ref{wallcrossingexamples}, Part 3]
    Let $\E\ra C\times\DD$ be a small family passing transversely through a generic stable bundle $E_0\in\W_g$, for which the total space of $\N_C^{g-1}(\E)\ra\DD$ is smooth by Proposition \ref{stablefamilysmoothness}. Proposition \ref{Wg_odp} shows that its central fibre has one ordinary double point. Picard--Lefschetz theory implies that the monodromy map $\Upsilon$ associated to the admissible loop $\E|_{\del\DD}$ is a generalised Dehn twist around the vanishing cycle, which is a (real) $(g-1)$-dimensional Lagrangian sphere. The formula \eqref{lefschetzdehntwist} is a consequence of equation \eqref{productcount} and the following general statement: Let $X$ be a closed symplectic manifold of dimension $2n$ with a Lagrangian sphere $L\subset X$. Then the generalised Dehn twist $\tau:X\ra X$ around $L$ has Lefschetz number
    \begin{equation}\label{symplecticlefschetz}
        \#\Fix(\tau)=\chi(X)-\chi(L).
    \end{equation}
    To prove this, consider the Picard--Lefschetz formula
    $$\tau_*(x)=x-(-1)^{n(n-1)/2}\langle x,[L]\rangle [L],$$
    where $x\in H_*(X;\Z)$ and $[L]\in H_n(X;\Z)$ is the fundamental class of $L$. This intersection pairing uses the symplectic orientation of $X$ (as in \cite[\S5a, (5.2)]{Seidel2000Graded}), for which we also have $\langle[L],[L]\rangle=(-1)^{n(n-1)/2}\chi(L)$. This sign corrects for the fact that the isomorphism $TL\simeq NL$ given by a compatible almost-complex structure is not necessarily orientation-preserving. Then $\tau_*[L]=(1-\chi(L))[L]$, so $\tr(\tau_*|_{H_n})=\dim H_n - \chi(L)$. The intersection pairing vanishes for classes of differing degrees, so $\tau_*|_{H_j}=\id$ for $j\neq n$ and formula \eqref{symplecticlefschetz} follows from the definition of the Lefschetz number.
\end{proof}

\section{Multi-vortices on genus 2 curves}\label{genus2}
Let us now specialise to the case where $C$ is a smooth projective curve of genus 2. It is a classical result of Narasimhan and Ramanan that the moduli space $\MM_C^{ss}$ can be identified with $\P^3$ \cite{NRgenus2}. The divisor $\W_1$ is the strictly semi-stable locus, which is a Kummer surface in $\P^3$ with 16 ordinary double points. The divisor $\W_2$ consists of 16 hyperplanes, classically known as tropes. These are arranged in the so-called (6,16) configuration, meaning that each trope contains 6 singularities of the Kummer, and each singularity is contained in 6 of the tropes. By Proposition \ref{wobblystructure}, the three interesting cases to consider are $d=0,1,2$, where
$$\W(2,0)=\W(2,2)=\W_1,\hspace{4mm}\W(2,1)=\W_1\cup\W_2.$$
The monodromy around these loci is described by Theorem \ref{genus2theorem}, which we will spend the rest of the paper proving. The computations required are more delicate than those in Section \ref{monodromysection}, since $\W_1$ consists of S-equivalence classes of strictly semi-stable bundles. We construct small families passing transversely through it, for which there are two choices of central fibre. This is because near $\W_1$, $\MM_C^{ss}$ is locally modelled by a GIT quotient.

\subsection{Generalised theta divisors and the (6,16) configuration}\label{616}

We begin by relating the moduli space $\N_C^1(E)$ to the generalised theta divisor. 

\begin{defn}
    The generalised theta divisor of a rank 2 bundle $E$ on a smooth genus $2$ curve $C$ is
    $$\theta(E)=\{L\in J^{1}_C:h^0(L\otimes E)>0\}.$$
\end{defn}

There is a surjective forgetful map $\pi:\N_C^{1}(E)\ra \theta(E)$ with projective space fibres. In this genus $2$ case, it can be shown that $\theta(E)$ is a actually a divisor in $J_C^1$\footnote{This fails for $E$ of higher rank and higher genus, as it may be the case that $\theta(E)=J^{g-1}_C$.}. It is linearly equivalent to twice the canonical theta divisor
$$\Theta=\{L\in J^{1}_C : h^0(L)>0\}.$$
The divisor $\theta(E)$ can be described geometrically as follows. The linear system $|2\Theta|$ is base-point free and defines a holomorphic map $\varphi:J^{1}_C\ra|2\Theta|^*$. The line bundle $\O(2\Theta)\ra J^{1}_C$ is isomorphic to $\varphi^*\O(H)$, where $H\subset |2\Theta|^*$ is a hyperplane. The point $\theta(E)\in|2\Theta|$ can be viewed as either the zero set of a section of $\O(2\Theta)$, or as the dual hyperplane $P_E\subset |2\Theta|^*$. This yields
$$\theta(E) = \varphi^{-1}(P_E).$$
Narasimhan and Ramanan proved the following beautiful relationship between generalised theta divisors and the geometry of $\MM_C^{ss}$.

\begin{theorem}[\cite{NRgenus2}]\label{nrgenus2}
    Let $C$ be a smooth genus $2$ curve. Then for every $E\in\MM_C^{ss}$, the generalised theta divisor $\theta(E)$ is linearly equivalent to $2\Theta$ and $\theta:\MM_C^{ss}\ra|2\Theta|$ is an isomorphism. The map $k:J^0_C\ra \MM_C^{ss}$ given by $k(A)=[A\oplus A^{-1}]$\footnote{We will use square brackets to denote the S-equivalence class of a strictly semi-stable bundle.} factors through the Kummer surface $J^0_C/\Z_2$, embedding it into $\MM_C^{ss}$. Its image is the strictly semi-stable locus $\W_1\subset \MM_C^{ss}$. 
\end{theorem}

The map $\varphi$ factors through the Kummer surface $J^1_C/\Z_2$, embedding it into $|2\Theta|^*$. The action of $\Z_2$ on $J_C^1$ is given by the involution $L\mapsto L^{-1}\otimes K_C$, so the singular points of $J^1_C/\Z_2$ are theta characteristics: bundles $L$ with $L^2=K$. The two Kummer surfaces $\W_1=k(J^0_C)\subset|2\Theta|$ and $\W^1:=\varphi(J_C^1)\subset|2\Theta|^*$ are related to each other in the following way. Every point $\varphi(L)\in\W^1$ defines a hyperplane in $|2\Theta|$, which via $\theta$ can be identified with
\begin{align*}
    H_L &= \{[E]\in\MM_C^{ss}:h^0(L\otimes E)>0\}\\
    &= \{\textnormal{non-trivial }E : 0\ra L^{-1}\ra E\ra L \ra 0\}.
\end{align*}
These hyperplanes are tangent to $\W_1$, and every tangent plane arises in this way. Dually, points of $\W_1$ correspond to tangent planes of $\W^1$. We will be particularly interested in $H_L$ when $L$ is one of the 16 theta characteristics. Classically, these hyperplanes are known as tropes, forming a $(6,16)$ configuration with the singularities of $\W_1$. This means that each of the 16 trope hyperplanes passes through 6 singular points of $\W_1$, and every one of the 16 singular points has 6 trope hyperplanes passing through it. This can be seen explicitly in $\MM_C^{ss}$ by noting that
\begin{equation}\label{kummertropeintersection}
    \W_1\cap H_L = \{[L(-p)\oplus L^{-1}(p)] : p\in C\}.
\end{equation}
There is a map $C\ra \W_1\cap H_L$ given by $p\mapsto [L(-p)\otimes L^{-1}(p)]$. Since $L^{2}\simeq K$, both $p$ and $\iota(p)$ map to isomorphic bundles, where $\iota:C\ra C$ is the hyperelliptic involution. We then see that $\W_1\cap H_L$ is isomorphic to the quotient of $C$ by $\iota$, which is a $\P^1$. This is a realisation of the classical fact that the intersection $\W_1\cap H_L$ is a double conic \cite{Dolgachev2012}. It contains exactly 6 of the 16 singular points of $\W_1$, which are the images of the Weierstrass points of $C$ under this map. To complete the picture, we briefly describe how two different tropes intersect. Suppose $L\neq L'$ are theta characteristics. Then
$$\W_1\cap H_L\cap H_{L'} = \{[L(-p)\oplus L^{-1}(p)], [L(-q)\oplus L^{-1}(q)]\},$$
where $p,q\in C$ are the unique points such that $L\otimes L'=\O(p+q)$ (uniqueness follows from $h^0(L\otimes L')=1$). Since $H_L\cap H_{L'}$ is a line, it is the unique line between these two points in $\W_1\cap H_L\cap H_{L'}$. 

\subsection{Local structure of $\MM_C^{ss}$ near the Kummer surface}

To study the monodromy of $\N_C$ induced by a loop $E_t$ passing around the Kummer surface $\W_1$, we need an analogue of Lemma \ref{stablefamily} for a disk passing transversally through a strictly semi-stable point. This is a more subtle problem, since such a disk does not determine a unique family of bundles and care needs to be taken when dealing with the central fibre. Every point in $\W_1$ is an S-equivalence class $[A\oplus A^{-1}]$ for some $A\in J_C^0$. Its singularities occur when $A^2=\O$, so a generic point has $A^2\neq \O$. The following computation describes the bundles in this S-equivalence class.

\begin{lemma}\label{sequivalence}
    Suppose $A^2\neq\O$. Then up to bundle isomorphism there are exactly three bundles in the S-equivalence class $[A\oplus A^{-1}]$; $A\oplus A^{-1}$ itself and two other bundles $E_\pm$, which are the unique non-trivial extensions
    \begin{equation}\label{Epm}
        0\ra A\xrightarrow{i} E_{+}\xrightarrow{j} A^{-1}\ra 0,\hspace{4mm}0\ra A^{-1}\xrightarrow{i} E_{-}\xrightarrow{j} A\ra 0.
    \end{equation}
    These bundles are simple and have $h^1(\End_0(E_{\pm}))=3$, while $A\oplus A^{-1}$ has
    $$\Aut(A\oplus A^{-1})=\C^*\times\C^*,\hspace{4mm} h^1(\End_0(A\oplus A^{-1}))=4.$$
\end{lemma}

\begin{proof}
    The bundles $E_{\pm}$ have Jordan--H\"{o}lder filtrations $0\subset A^{\pm1}\subset E_\pm$, with graded object $\gr\simeq A\oplus A^{-1}$. This shows they are S-equivalent to $A\oplus A^{-1}$. Conversely, suppose $E$ is S-equivalent to $A\oplus A^{-1}$. Then a Jordan--H\"{o}lder filtration $0\subset E_1\subset E$ has $\gr=E_1\oplus E/E_1 \simeq A\oplus A^{-1}$, which identifies $E_1$ with either $A$ or $A^{-1}$. Then $E$ is an extension of $A$ by $A^{-1}$ (or vice versa), and is isomorphic to $A\oplus A^{-1}$ or $E_{\pm}$. Now consider the short exact sequences
    \begin{equation}\label{shortexact1}
        0\ra \End_0^A(E_+)\ra \End_0(E_+)\xrightarrow{p} \Hom(A,A^{-1})\ra 0,
    \end{equation}
    \begin{equation}\label{shortexact2}
        0\ra A^2\ra \End_0^A(E_+)\ra \O\ra 0.
    \end{equation}
    Here, $\End_0$ denotes trace-free endomorphisms and $\End_0^A$ denotes trace-free endomorphisms preserving the sub-bundle $A$. The map $p$ sends $\theta:E_+\ra E_+$ to the composition $j\circ\theta\circ i:A\ra A^{-1}$. Taking cohomology and using $h^0(A^{\pm2})=0$, we obtain
    \begin{equation}\label{les1}
        0\ra H^0(\End_0^A(E_+))\ra H^0(\End_0(E_+))\ra 0,
    \end{equation}
    \begin{equation}\label{les2}
        0\ra H^1(\End_0^A(E_+))\ra H^1(\End_0(E_+))\xrightarrow{p} H^1(A^{-2})\ra 0,
    \end{equation}
    \begin{equation}\label{les3}
        0\ra H^0(\End_0^A(E_+))\ra H^0(\O)\xrightarrow{e}H^1(A^2)\ra H^1(\End_0^A(E_+))\ra H^1(\O)\ra 0.
    \end{equation}
    The coboundary map in \eqref{les3} sends the constant section $1$ to the non-zero extension class $e\in H^1(A^{2})$ of $E_+$, so \eqref{les1} shows that $h^0(\End_0(E_+))=h^0(\End_0^A(E_+))=0$ and $E_+$ is simple. It follows from Riemann--Roch that $h^1(A^{\pm2})=1$, so \eqref{les2} and \eqref{les3} imply $h^1(\End_0(E_+))=3$. The same results hold for $E_-$. As for $A\oplus A^{-1}$, $h^1(A^{\pm2})=1$ directly implies $h^1(\End_0(A\oplus A^{-1}))=4$. Similarly, one can show that $H^0(\End(A\oplus A^{-1}))=H^0(\O)\oplus H^0(\O)$ is two-dimensional. Requiring invertability shows that $\Aut(A\oplus A^{-1})=\C^*\times\C^*$.
\end{proof}

\begin{prop}\label{semistablefamily}
    Suppose $u:\DD\ra\MM_C^{ss}$ is a holomorphic disk passing transversely through a point $[A\oplus A^{-1}]\in\W_1$ with $A^2\neq \O$, such that $u(t)$ is stable for $t\neq 0$. Then up to isomorphism, there are exactly two families $\E\ra C\times\DD$ such that $\E|_{C\times\{t\}}\simeq u(t)$ for $t\neq 0$ and $\E|_{C\times\{0\}}$ is semi-stable. The semi-stable central fibres of these two families $\E_+$ and $\E_-$ are $E_{+}$ and $E_{-}$, as defined in \eqref{Epm}. Moreover, the Kodaira--Spencer maps $\kappa_{\pm}:T_0\DD\ra H^1(\End_0(E_\pm))$ of $\E_\pm$ satisfy $\im(\kappa_{\pm})\cap\ker(p)=0$, where $p$ is the map in the exact sequence \eqref{les2}. 
\end{prop}

\begin{proof}
    We first prove the existence of $\E_+$ (the same argument works for $\E_-$). The vector bundle $E_+$ can be described using a Dolbeault operator $\delbar_+:\Omega^{0}(\underline{\C}^2)\ra\Omega^{0,1}(\underline{\C}^2)$ on the trivial vector bundle $\underline{\C}^2\ra C$. Let $H$ be the subspace of $\Omega^{0,1}(\End_0(\underline{\C}^2))$ consisting of $\delbar_{+}$-harmonic forms. The inclusion map $h:H\ra \Omega^{0,1}(\End_0(\underline{\C}^2))$ determines a family of Dolbeault operators $\delbar_{s} = \delbar_+ + h(s)$ on $\underline{\C}^2$. This determines a family of holomorphic vector bundles $\F\ra C\times H$ with $\F_0=E_+$. Fix an identification of $H$ with $H^1(\End_0(E_+))$, so that the Kodaira--Spencer map $\kappa_{\F}$ of $\F$ at $0$ is the identity. 

    Since semi-stability is an open condition, there is a neighbourhood $S\subset H$ of $0$ such that the map $f:S\ra \MM_C^{ss}$ given by $f(s)=[\F_s]$ is well-defined. The map assigning $s\in S$ to the isomorphism class of $\F_s$, rather than the S-equivalence class, is a local homeomorphism \cite[Thm.~7.3.17]{kobayashi}. To show that $f$ itself is injective, suppose $\F_s$ and $\F_{s'}$ are S-equivalent. If they are stable, then they are isomorphic and $s=s'$, so suppose they are strictly semi-stable and S-equivalent to $B\oplus B^{-1}$. Since $\Hom(A^{-1},E_+)=0$, Grauert's semicontinuity theorem \cite[Thm.~III.12.8]{hartshorne} implies that $\Hom(B^{-1},\F_{s})=0$ for $(s,B)$ sufficiently close to $(0,A)$. Lemma \ref{sequivalence} and the strict semi-stability of $\F_s$ show that it is the unique non-trivial extension of $B^{-1}$ by $B$. So after possibly shrinking $S$, $\F_s\simeq\F_{s'}$ and $s=s'$. Then $f$ is an injective holomorphic map, so it is a biholomorphism onto its image.

    For $W=f^{-1}(\W_1)\subset S$, we will now show that $T_0 W=\ker(p)$. Let $s:\DD\ra W$ be a holomorphic function with $s(0)=0$, for which $s'(0)\in T_0 W$. Up to isomorphism, this determines a family of line bundles $\mathcal{A}\ra C\times\DD$ with $\A_0=A$. The same argument as above shows that $\F_{s(\tau)}$ is the non-trivial extension of $\A_\tau^{-1}$ by $\A_\tau$. This has $\dim\Hom(\A_\tau,\F_{s(\tau)})=1$ for all $\tau$, so $\Hom(\A,s^*\F)\ra\DD$ is a line bundle by \cite[Thm.~III.12.11]{hartshorne}. The inclusion $i:A\hookrightarrow E_+$ extends to a non-zero holomorphic section $\varphi:\A\ra s^*\F$. At $\tau=0$, $\varphi$ is compatible with the Kodaira--Spencer maps of the families, so we have
    $$p(s'(0))=p(\kappa_{\F}(s'(0)))=j\circ\kappa_{\F}(s'(0))\circ i = j\circ i\circ \kappa_{\A}(\del_\tau)=0.$$
    The inclusion $T_0 W\subset\ker(p)$ is an equality since both spaces are 2-dimensional. Let $\tilde{u}=f^{-1}\circ u$ and define $\E_+ = (1\times\tilde{u})^*\F$. The Kodaira--Spencer map of this family is $\kappa_{+}=d\tilde{u}$, and since $\tilde{u}$ is transverse to $W$, we have
    \begin{equation}\label{transversalitycondition}
        \im(\kappa_+)\cap\ker(p)=0.
    \end{equation}
    The uniqueness up to isomorphism of $\E_\pm$ follows from the fibres being simple. To see that there is no family agreeing with $u(t)$ for $t\neq 0$ and with central fibre $A\oplus A^{-1}$, suppose $\E\ra C\times\DD$ is such a family. Since $\Hom(A\oplus A^{-1},E_+)$ is spanned by $i\circ\pr_A$, we have $\dim\Hom(\E_t,(\E_+)_t)=1$ for all $t$. Then as before, there is a holomorphic map $\varphi:\E\ra\E_+$ with $\varphi_0=i\circ\pr_A$. However,
    $$p(\kappa_+(\del_t))\circ\pr_A = j\circ\kappa_+\circ\varphi_0=j\circ\varphi_0\circ\kappa_{\E}(\del_t)=0,$$
    which contradicts \eqref{transversalitycondition}. 
\end{proof}

\begin{remark}
    Luna's slice theorem \cite{Luna1973} (see also \cite{Laszlo1996}) offers another way to prove Proposition \ref{semistablefamily}. It identifies a neighbourhood of $[A\oplus A^{-1}]$ in $\M_C^{ss}$ with the GIT quotient $V\git \C^*$, where $V=H^1(A^2)\oplus H^1(A^{-2})\oplus H^1(\O)$ has coordinates $(x,y,z)$ and $\C^*$ acts with weights $(1,-1,0)$. If $(w=xy,z)$ are coordinates on $V\git\C^*\simeq \C^3$, a holomorphic disk intersecting $\W_1=\{w=0\}$ transversely at $0$ would look like $t\mapsto(t,0)$. There are two lifts $(1,t,0)$ or $(t,1,0)$ to $V$, which correspond to the families $\E_\pm$. The central fibre $A\oplus A^{-1}$ can only occur for a disk like $t\mapsto (t^2,0)$ because it has a lift $(t,t,0)$ vanishing at $t=0$. However, such a disk does not intersect $\{w=0\}$ transversely.
\end{remark}

\subsection{Monodromy around $\W(2,1)$ for $d=1$}

For a genus $2$ curve, the most interesting situation to study is when $d=g-1=1$ and $\W(2,1)=\W_1\cup\W_2$ is the full wobbly divisor. We find non-trivial monodromy around both $\W_1$ and $\W_2$. The proof of the next lemma can be found in \cite{NRgenus2}.

\begin{lemma}\label{genus2lemma}
    Suppose $C$ is a smooth genus $2$ curve and let $E$ be a stable rank $2$ bundle with trivial determinant. Then for every $L\in J_C^1$,
    \begin{enumerate}
        \item $h^0(L\otimes E)\leq 1$.
        \item Any non-zero map $L^{-1}\ra E$ is everywhere injective. 
    \end{enumerate}
\end{lemma}

The relationship between $\N$ and the generalised theta divisor can be used to understand its geometry when $E\notin\W(2,1)$. 

\begin{prop}[\cite{Doan_2018}]\label{genus5}
    Suppose $E\notin\W(2,1)$. Then both $\N_C^1(E)$ and $\M_C^1(E)$ are biholomorphic to $\theta(E)$, which is a smooth genus $5$ curve.
\end{prop}

\begin{proof}
    The map $\M_C^1(E)\ra \theta(E)$ is an isomorphism (and so is $\pi$) if and only if the fibre over any $L\in\theta(E)$ consists of just one point. This is equivalent to showing that $h^0(E\otimes L)=1$ and for non-zero $\alpha\in H^0(E\otimes L)$, there is no non-zero $\beta\in H^0(E^*\otimes L^{-1}\otimes K)$ such that $\beta\alpha=0\in H^0(K)$. Now suppose $E\notin \W(2,1)$, so it is stable and Lemma \ref{genus2lemma} guarantees that $h^0(E\otimes L)=1$ for all $L\in\theta(E)$. Suppose there is $\beta\neq 0$ with $\beta\alpha=0$. Also by Lemma \ref{genus2lemma}, $\alpha:L^{-1}\ra E$ is injective and $\beta:E\ra L^{-1}\otimes K$ is surjective, so we have the exact sequence
    $$0\ra L^{-1}\xrightarrow{\alpha}E\xrightarrow{\beta}L^{-1}\otimes K\ra 0.$$
    The condition $\det(E)\simeq\O$ requires $L^2\simeq K$, but $E\in H_L$ contradicts $E\notin\W(2,1)$. Therefore,
    $$\M_C^1(E)\simeq \N_C^1(E)\simeq \theta(E).$$
    By our previous remarks, we have
    $$\theta(E) = \varphi^{-1}(P_E) = k^{-1}(\W^1\cap P_E).$$
    For any $E\in H_L$, the dual hyperplane $P_E$ passes through the singular point $\varphi(L)\in \W^1$. If $E\notin \W(2,1)$, $P_E$ does not contain any of the singular points of $\W^1\subset|2\Theta|^*$ and is also not tangent to $\W^1$. Then $P_E$ intersects $\W^1$ transversally and $\W^1\cap P_E$ is a smooth quartic plane curve of genus 3. Its preimage under $k$ is the theta divisor $\theta(\E)$, which is a smooth genus 5 curve by the Riemann-Hurwicz formula. 
\end{proof}

We now turn to understanding $\N_C^1(E)$ when $E\in\W_2$ is chosen generically. Indeed, Proposition \ref{Wg_odp} shows that for a generic stable $E\in\W_2$, $\N_C^1(E)$ has a unique ordinary double point. The next result uses the same ideas as Proposition \ref{genus5} to describe its geometry more explicitly.

\begin{prop}\label{tropefibre}
    Suppose $E\in H_L\subset \W_2$ is chosen generically, where $L^2=K$. Then $\N_C^1(E)$ is biholomorphic to $\theta(E)$, which is a curve of arithmetic genus $5$ with a node at $L$.
\end{prop}

\begin{proof}
    Suppose that
    $$E\in H_L\setminus\Big(\W_1\cup\bigcup_{M^2=K,M\neq L}H_M\Big),$$
    which is dense in $H_L$. Then since $E$ is stable, the argument above shows that $\theta(E)\simeq\N_C^1(E)$. By construction, the only singular point of $\W^1$ contained in the dual plane $P_{E}\subset|2\Theta|^*$ is $\varphi(L)$, and $P_{E}$ intersects $\W^1$ transversally elsewhere because $[E]\notin\W_1$. The intersection $X=\W^1\cap P_{E}$ is then a degree 4 plane curve with a unique node at $L$, so it has arithmetic genus $3$ and geometric genus $3-1=2$. We wish to determine $\theta(E)=k^{-1}(X)$, where $k:J^1_C\ra\W^1$ is branched double cover. The double cover $\theta(\E)\ra X$ branches only at the node $\varphi(L)$. Consider the normalisation map $\nu:Y\ra X$ which maps two distinct points $p,q$ on a smooth genus 2 curve $Y$ to $\varphi(L)\in X$, but is otherwise an isomorphism. There is a unique double cover $Z\ra Y$ branching at $p$ and $q$, which is a smooth genus 4 curve by the Riemann--Hurwicz formula. The following diagram commutes:
        \[\begin{tikzcd}
	    Z & {\theta(\E)} \\
	    Y & X
	    \arrow["{\nu'}", from=1-1, to=1-2]
	    \arrow[from=1-1, to=2-1]
	    \arrow[from=1-2, to=2-2]
	    \arrow["\nu", from=2-1, to=2-2]
        \end{tikzcd}\]
    Here, $\nu'$ is the normalisation map gluing the preimages of $p$ and $q$ in $Z$ together. This means that $\theta(\E)$ is a curve with arithmetic genus 5 and one node at $L$.
\end{proof}

The next case to consider is when $[E]$ is chosen generically in $\W_1$. The geometry of $\N_C^1(E)$ is more subtle in this case, so we need a series of preparatory lemmas before stating our main results in Lemma \ref{kummerfibre} and Proposition \ref{W1smoothness}. 

\begin{lemma}\label{generickummer}
    Suppose $A\in J_C^0$ and let $\iota:C\ra C$ be the hyperelliptic involution. Then $[A\oplus A^{-1}]\in \W_1\setminus(\W_1\cup\W_2)$ if and only if $A^2\simeq\O(x-y)$ for $x,y\in C$ with $y\notin\{x,\iota(x)\}$. 
\end{lemma}

\begin{proof}
    Recall that on the genus 2 curve $C$, a degree 2 divisor $D$ has $\dim|D|=1$ if and only if $D$ is canonical, otherwise $\dim|D|=0$. We claim that for $x,y,a,b\in C$, the following holds:
    \begin{equation}\label{linearequiv}
        x-y\sim a-b\textnormal{ if and only if }(x,y)\in\{(a,b),(\iota(b),\iota(a))\}.
    \end{equation}
    Certainly if $(x,y)=(\iota(b),\iota(a))$ we have the linear equivalence, since then $x+b$ and $y+a$ are both canonical divisors. On the other hand, if $x+b$ is not canonical then $\dim|x+b|=0$ forces $x+b=y+a$ and therefore $(x,y)=(a,b)$. 

    Now consider the map $\rho:C\times C\ra J_C^0$ given by $\rho(x,y)=\O(x-y)$. Its image is closed, so the open mapping theorem implies it is surjective and we can write $A^2\simeq\O(x-y)$ for some $x,y\in C$. We proceed by contrapositive. Suppose that $[A\oplus A^{-1}]\in\W_1\cap \W_2$, so by \eqref{kummertropeintersection} there exist $L\in J_C^1$ and $p\in C$ with $L^2\simeq K$ and $A\simeq L(-p)$. Then $\O(x-y)\simeq K(-2p)\simeq\O(\iota(p)-p)$ and \eqref{linearequiv} implies that $x=\iota(p)$ and $y=p$. Conversely, if $y=x$ we have $[A\oplus A^{-1}]\in \W_1\cap\W_2$ because $A^2=\O$. If $y=\iota(x)$, let $L=A(i(x))$ which satisfies $L^2\simeq \O(x+\iota(x))\simeq K$. Then \eqref{kummertropeintersection} shows that $[A\oplus A^{-1}]\in \W_1\cap H_L\subset \W_1\cap\W_2$.
\end{proof}

\begin{lemma}\label{thetatangentspace}
    Identify every tangent space of $J_C^1$ with $H^1(\O)$ and let $\phi_K:C\ra|K|$ be the canonical map. Then Serre duality identifies $T_{\O(p)}\Theta\subset H^1(\O)$ with the line $\phi_K(p)\subset H^0(K)^*$.
\end{lemma}

\begin{proof}
    Let $u:C\ra J_C^1$ be the Abel--Jacobi map $u(p)=\O(p)$, whose image is $\Theta$. Consider the short exact sequence
    $$0\ra \O\ra \O(p)\ra \O(p)|_p\ra 0.$$
    The tangent space $T_pC$ can be canonically identified with $H^0(\O(p)|_p)$ (see \cite{acgh}) and the derivative $d_pu:T_pC\ra H^1(\O)$ is identified with the coboundary map $\delta:H^0(\O(p)|_{p})\ra H^1(\O)$. This implies that $T_{\O(p)}\Theta=\im(\delta)$, which is the kernel of $H^1(\O)\ra H^1(\O(p))$. This map is Serre dual to the inclusion $H^0(K(-p))\ra H^0(K)$, so its kernel can be identified with the annihilator
    $$H^0(K(-p))^0\subset H^0(K)^*\simeq H^1(\O).$$ 
    To relate this to the canonical map, consider the evaluation map $\ev_p:H^0(K)\ra K_p$. Its dual globalises to $\ev^*:K^{-1}\ra H^0(K)$ and the canonical map is $\phi_K(p)=[\ev_p^*]$. Explicitly, it is given by $\ev_p^*:\xi\mapsto (s\mapsto \langle s(p),\xi\rangle)$, so the line $[\ev_p^*]$ consists of functionals vanishing on sections of $H^0(K)$ which are zero at $p$. This is precisely the annihilator $H^0(K(-p))^0$, which identifies $T_{\O(p)}\Theta$ with $\phi_K(p)$.
\end{proof}

\begin{lemma}\label{kummerfibre}
    Suppose $[A\oplus A^{-1}]\in\W_1\setminus (\W_1\cup\W_2)$ and let $E_\pm$ be the two non-trivial extensions \eqref{Epm} determined by $A\in J_C^0$. Then $\N_C^1(E_{\pm})$ is isomorphic to $\theta(E_{\pm})$, which is a curve of arithmetic genus $5$ with two nodes. 
\end{lemma}

\begin{proof}
    By Lemma \ref{generickummer}, there are $x,y\in C$ with $y\notin\{x,\iota(x)\}$ such that $A^2\simeq \O(x-y)$. Without loss of generality, we work with $E_+$. Let $\Theta_{A^{\pm1}}=\{A^{\pm1}(p) : p\in C\}$ be the translation of the canonical theta divisor $\Theta$ by $A^{\pm1}$. We first prove that $\theta(E_+)=\Theta_A + \Theta_{A^{-1}}$. Suppose $L\in\theta(E_+)$ and consider the long exact sequence
    \begin{equation}\label{thetalongexact}
        0\ra H^0(L\otimes A)\ra H^0(L\otimes E_+)\ra H^0(L\otimes A^{-1})\ra H^1(L\otimes A)\ra \cdots
    \end{equation}
    Since $h^0(L\otimes E_+)\geq 1$, $h^0(L\otimes A)$ and $h^0(L\otimes A^{-1})$ cannot both be zero. Then $L=A^{\pm1}(p)\in\Theta_{A^{\pm1}}$ for some $p\in C$. Conversely, if $L=A^{-1}(p)$ for some $p\in C$, the map $H^0(p)\ra H^0(L\otimes E_+)$ is injective and $L\in\theta(E_+)$. On the other hand, if $L=A(p)$, the exact sequence \eqref{thetalongexact} becomes
    $$0\ra H^0(x-y+p)\ra H^0(L\otimes E_+)\ra H^0(p)\ra H^1(x-y+p)\ra \cdots$$
    The first term is non-zero if $p\in\{y,\iota(x)\}$, which ensures that $L\in\theta(E_+)$, but it vanishes when $p\notin\{y,\iota(x)\}$. In that case, Riemann--Roch guarantees that $h^1(x-y+p)=0$, so $h^0(L\otimes E_+)=1$ and $L\in\Theta(E_+)$. 
    
    It follows from \eqref{linearequiv} that $\Theta_A$ meets $\Theta_{A^{-1}}$ only at the two points
    $$A(y)\simeq A^{-1}(x),\hspace{4mm}A(\iota(x))\simeq A^{-1}(\iota(y)).$$
    To prove that these points are nodes, it suffices to check that $\Theta_A$ and $\Theta_{A^{-1}}$ intersect transversally there. Since translation by $A^{\pm1}$ acts as the identity on $H^1(\O)$, Lemma \ref{thetatangentspace} implies $T_{A^{\pm1}(p)}\Theta_{A^{\pm1}}=\varphi_K(p)$. It follows that $T_{A(y)}\Theta_A=T_{A^{-1}(x)}\Theta_{A^{-1}}$ if and only if $\varphi_K(x)=\varphi_K(y)$, which is equivalent to $y\in\{x,\iota(x)\}$. This is false by assumption, so in fact $T_{A(y)}\Theta_A\neq T_{A^{-1}(x)}\Theta_{A^{-1}}$ are transverse and $A(y)\in\theta(E_+)$ is a node (the same argument works for the other point).

    It remains to prove that $h^0(L\otimes E_+)=1$ for every $L\in\theta(E_+)$, which will show that $\N_C^1(E_+)$ is biholomorphic to $\theta(E_+)$ (Lemma \ref{genus2lemma} does not apply because $E_+$ is not stable). We already showed that $h^0(L\otimes E_+)=1$ when $L=A(p)$ and $p\notin\{y,\iota(x)\}$, and a similar argument works when $L=A^{-1}(p)$ and $p\notin\{x,\iota(y)\}$. Now let $L=A(y)$, so that the exact sequence \eqref{thetalongexact} is
    $$0\ra H^0(x)\ra H^0(L\otimes E_+)\ra H^0(y)\ra H^1(x)\ra \cdots$$
    To prove $h^0(L\otimes E_+)=1$, it suffices to show that the coboundary map $H^0(y)\mapsto H^1(x)$ is an isomorphism. It is the cup product with a non-zero section in $H^1(x-y)$ (all three of these spaces are 1-dimensional). Consider the exact sequence
    $$0\ra \O(x-y)\ra \O(x)\ra \O(x)|_y\ra 0.$$
    The induced map $H^1(x-y)\ra H^1(x)$ is the cup product with a non-zero section in $H^0(y)$, and its cokernel is $H^1(\O(x)|_y)=0$. Then it is an isomorphism, and so is the original coboundary map. The same argument works for $L=A(\iota(x))$, so $\N_C^1(E_+)\simeq\theta(E_+)$. We conclude by remarking that swapping $E_+$ for $E_-$ interchanges $\Theta_{A}$ and $\Theta_{A^{-1}}$, so $\theta(E_+)\simeq\theta(E_-)$.
\end{proof}

\begin{remark}
    Recall that $E_{+}$ and $E_{-}$ are S-equivalent to $A\oplus A^{-1}$ (c.f. Lemma \ref{sequivalence}). It is not hard to show that the generalised theta divisor $\theta(A\oplus A^{-1})$ is biholomorphic to $\theta(E_{\pm})$, however the same is not true for $\N_C^{1}(A\oplus A^{-1})$. The map $\pi:\N_C^{1}(A\oplus A^{-1})\ra \theta(A\oplus A^{-1})$ has point fibres everywhere except the two nodes, over which the fibre is $\P^1$. By analysing how they glue, one finds that $\N_C^1(A\oplus A^{-1})$ consists of two copies of $C$ (that is, $\Theta_{A^{\pm1}}$) with two copies of $\P^1$ spanning from $A(y)\in\Theta_A$ to $A^{-1}(x)\in\Theta_{A^{-1}}$ and from $A(\iota(x))\in\Theta_{A}$ to $A^{-1}(\iota(y))\in\Theta_{A^{-1}}$. 
\end{remark}

Theorem \ref{nollauthm} does not apply to $\W_1$ for genus 2 curves, however the next lemma shows that the same result still holds in this case. In analogy with Proposition \ref{stablefamilysmoothness}, we use it to prove Proposition \ref{W1smoothness} below. This establishes smoothness of $\N_C^1(\E_\pm)$ for the semi-stable families $\E_\pm\ra C\times\DD$ (see Proposition \ref{semistablefamily}), which is essential for applying Picard--Lefschetz theory to understand the monodromy. 

\begin{lemma}\label{uniquehiggsW1}
    Let $E_\pm$ be the two non-trivial extensions \eqref{Epm} determined by $A\in J_C^0$ with $A^2\neq\O$. Then there is a unique non-zero nilpotent Higgs field for each of $E_\pm$.
\end{lemma}

\begin{proof}
    Without loss of generality, we work with $E_+$. By Remark \ref{wobblyequivalence}, every non-zero nilpotent Higgs field for $E_+$ arises as a composition 
    \begin{equation}\label{nilpotenthiggs}
        E_+\xrightarrow{j}M^{-1}\xrightarrow{\sigma}M\otimes K\xrightarrow{i\otimes 1}E_+\otimes K,
    \end{equation}
    where $M$ is a degree $0$ line sub-bundle of $E_+$ admitting a non-zero map $\sigma\in\Hom(M^{-1},M\otimes K)$. Tensoring \eqref{Epm} by $M^{-1}$, we get
    $$0\ra H^0(M^{-1}\otimes A)\ra H^0(M^{-1}\otimes E_+)\ra H^0(M^{-1}\otimes A^{-1})\ra H^1(M^{-1}\otimes A)\ra\cdots$$
    For $h^0(M^{-1}\otimes E_+)$ to be non-zero, we require either $h^0(M^{-1}\otimes A)>0$ or $h^0(M^{-1}\otimes A^{-1})>0$. In the first case $M\simeq A$, and in the second case $M\simeq A^{-1}$ and $H^0(A\otimes E_+)\ra H^0(\O)$ is surjective. However, that means that the coboundary map $H^0(\O)\ra H^1(A^2)$ is zero, which is not the case because it sends the constant section $1$ to the non-zero extension class of $E_+$. Thus the only choice is $M\simeq A$, in which case $h^0(A^2\otimes K)=1$ because $A^2\neq \O$. This means that the non-zero nilpotent Higgs field $\phi:E_+\ra E_+\otimes K$ is unique up to scaling.
\end{proof}

\begin{prop}\label{W1smoothness}
    Let $\E_{\pm}\ra C\times\DD$ be the two families determined by a holomorphic disk transversally intersecting $\W_1$ at $[A\oplus A^{-1}]\in\W_1\setminus(\W_1\cap\W_2)$, as described by Proposition \ref{semistablefamily}. Then the total space of $\N_C^1(\E_{\pm})\ra \DD$ is smooth.
\end{prop} 

\begin{proof}
    Without loss of generality, we work with $\E_+$. The argument of Proposition \ref{stablefamilysmoothness} shows that $\N_C^1(\E_+)$ is smooth at any point $(L,\alpha,t)$ such that $(L,\alpha)$ is a smooth point of $\N_C^{1}(E_t)$. For each of the two singular points $(L_1=A(y),\alpha_1)$, $(L_2=A(\iota(x)),\alpha_2)$ of $\N_C^1(E_+)$, there is $\beta_i\in H^0(L_i^{-1}\otimes E_+\otimes K)$ with $\beta_i\alpha_i=0$. Lemma \ref{uniquehiggsW1} shows that there is a unique non-zero nilpotent Higgs field $\phi\in H^0(\End_0(E_+)\otimes K)$, which must agree with $\alpha_i\beta_i$ up to rescaling. Then $\mu$ \eqref{mu} is given by
    $$\mu(\beta_i)=(\kappa^*(\phi),0).$$
    We need to show that $\mu$ is injective, so it suffices to check that the Kodaira--Spencer map $\kappa:T_0\DD\ra H^1(\End_0(E_+))$ for the family $\E_+$ satisfies $\kappa^*(\phi)\neq 0$. From the sequence \eqref{les2}, consider the map $p:H^1(\End_0(E_+))\ra H^1(A^{-2})$. Its adjoint is
    $$p^*:H^0(A^2\otimes K)\ra H^0(\End_0(E_+)\otimes K),$$
    which is interpreted via Serre duality. Explicitly, $p^*$ is given by the composition \eqref{nilpotenthiggs}, so its 1-dimensional image is spanned by $\phi$ and we have
    \begin{equation}\label{complement}
        \ker(p)=\ker(\phi^*).
    \end{equation}
    Proposition \ref{semistablefamily} shows that $\im(\kappa_+)\cap\ker(p)=0$, so $\phi^*(\kappa_+(\del_t))\neq 0$. This is equivalent to $\kappa^*(\phi)\neq 0$, so $\mu$ is injective and $\N_C^1(\E_+)$ is smooth. 
\end{proof}

We can now combine the results of this section to prove parts 2a and 2b of Theorem \ref{genus2theorem}. This completely describes the monodromy of $\N_C^1$ around small loops linking $\W(2,1)$. 

\begin{proof}[Proof of parts 2a and 2b of Theorem \ref{genus2theorem}]
    That $\N_C^1(E)$ is a smooth genus 5 curve for $E\notin\W(2,1)$ is the content of Proposition \ref{genus5}. For $i\in\{1,2\}$, let $u_i:\DD\ra \MM_C^{ss}$ be holomorphic disks passing transversely through generic points $u_i(0)\in\W_i$, with $u(t)\notin\W(2,1)$ for all $t\neq 0$. Lemma \ref{stablefamily} and Proposition \ref{semistablefamily} determine families of bundles $\E_i\ra C\times\DD$ such that $[\E_i|_{C\times\{t\}}]= u_i(t)$ for all $t$. The associated families of multi-vortex moduli spaces $\N_C^1(\E_i)\ra\DD$ are smooth by Propositions \ref{stablefamilysmoothness} and \ref{W1smoothness}. The central fibre of $\N_C^1(\E_2)$ has one node by Proposition \ref{tropefibre}, so Picard--Lefschetz theory implies that the monodromy map $\Upsilon_2$ is a Dehn twist around vanishing cycle corresponding to this node. Proposition \ref{kummerfibre} shows that the central fibre of $\N_C^1(\E_1)$ consists of two copies of $C$ (that is, $\Theta_A$ and $\Theta_{A^{-1}}$) joined by two nodes. The vanishing cycles of these distinct nodes are disjoint, so the monodromy map $\Upsilon_1$ is a composition of the (commuting) Dehn twists around each of them.
\end{proof}

\begin{remark}
    Upon taking cohomology, parts 2a and 2b Theorem \ref{genus2theorem} define a rank 10 local system on $\MM_C^{ss}\setminus\W(2,1)$. The geometric Langlands correspondence studies such local systems and their relationship to local systems on $C$ itself. Although there is no reason for them to be related, we can compare our local system to the rank 8 local system $\mathcal{F}_{dR,0}$ constructed in \cite[p.~17]{dps}. Its monodromy around each hyperplane in $\W_2$ is a transvection, and around $\W_1$ the monodromy is a composition of four commuting transvections. Dehn twists are transvections, so the two local systems share a vaguely similar structure. 
\end{remark}

\subsection{Monodromy around $\W_1$ for $d=0$ and $d=2$}

In the previous section, we found that for $d=1$, the monodromy around $\W(2,1)$ is extremely rich. In this section, we prove part 1 of Theorem \ref{genus2theorem}, which shows that for $d=0$ and $d=2$, the fixed point count for admissible loops passing around $\W(2,d)=\W_1$ is zero. When $d=0$, it is easy to see that $\N_C^0(E)$ is empty whenever $E\notin\W_1$ (that is, $E$ is stable). If not, there would be a point $(L,\alpha)\in\N_C^0(E)$, but $\alpha$ realises $L^{-1}$ as a destabilising sub-bundle. This shows that there are no fixed points in this case, for any admissible loop. 

Recall that as in ordinary Seiberg--Witten theory, there is a duality arising from conjugation of the spin$^c$ structure:
$$\#\M_Y(\mathfrak{s},\mathfrak{p})=\pm\#\M_Y(\overline{\mathfrak{s}},\mathfrak{p}).$$
The conjugate of the spin$^c$ structure $(d,\gamma)$ on $Y=S^1\times C$ is $(2g-2-d,-\gamma)$. Then Conjecture \ref{bijection} suggests that for $g=2$ and $d=2$, there are also no fixed points for any admissible loop. From the point of view of complex geometry, this is not immediately obvious. We prove it manually, by first showing that the moduli space $\N_C^2(E_\pm)$ has one ordinary double point. The monodromy is then a Dehn twist, and we show that the Lefschetz number of any composition of such Dehn twists vanishes. This, along with an explicit understanding of the geometry of $\N_C^2$, proves that the fixed point counts are zero for any admissible loop. 

\begin{lemma}\label{dequalstwo}
    Suppose $A\in J^0_C$ has $A^2\neq\O$ and $E_{\pm}$ are the two extensions given by \eqref{Epm}. Then $\N_C^{2}(E_{\pm})\ra J_C^{2}$ is a smooth $\P^{1}$-bundle over $J_C^{2}\setminus\{A^{\pm1}\otimes K\}$. It contains one ordinary double point in the exceptional fibre $\P H^0(A^{\pm1}\otimes K\otimes E_{\pm})=\P^{2}$ over $L=A^{\pm1}\otimes K$. 
\end{lemma}

\begin{proof}
    Without loss of generality, we work with $E_+$. Computations similar to those in the proof of Lemma \ref{kummerfibre} yield
    $$h^0(L\otimes E_+)=\begin{cases}
        3 & \textnormal{if }L\simeq A\otimes K,\\
        2 & \textnormal{otherwise. }  
    \end{cases}$$
    This shows that $\N_C^2(E_+)\ra J_C^2$ is a smooth $\P^1$-bundle away from $A\otimes K$, whose fibre is $\P^2$. Any singularity $(L,\alpha)$ must then occur when $L=A\otimes K$. Consider the exact sequence
    \begin{equation}\label{d2ses}
        0\ra H^0(A^2\otimes K)\ra H^0(A\otimes K\otimes E_+)\ra H^0(K)\ra 0,
    \end{equation}
    where as in the proof of Lemma \ref{uniquehiggsW1}, $H^0(A^2\otimes K)=\Hom(L^{-1},A)$ is spanned by $\sigma$. Its image in $H^0(A\otimes K\otimes E_+)=\Hom(L^{-1},E_+)$ is $\alpha=i\circ\sigma$, which has $\beta\alpha=0$ for $$\beta=j\in H^0(L^*\otimes E_+^*\otimes K)=\Hom(E_+,A).$$
    This shows that $\N_C^2(E_+)$ is singular at $(L,\alpha)=(A\otimes K,i\circ\sigma)$. Since $\Hom(E_+,A)$ is spanned by $j$ and $\ker(j\cdot)=\C\alpha$, there are no other singular points. To prove that it is an ordinary double point, we will use Proposition \ref{ordinarydoublepoint}. In our present situation, the line bundle $M$ appearing in the proof of Theorem \ref{wobblystructure} is $A$, and the divisors $D_i$ are $D_1=x+\iota(y)$ (where $A^2=\O(x-y)$) and $D_2=0$. Consider the following diagram:

    \[\begin{tikzcd}[column sep=small]
        && 0 &&&& \\
        && {H^0(L|_{D_1})} \\
        0 & {H^0(L\otimes E_+)/\mathbb{C}\alpha} & {H^0(L\otimes Q)} & {H^1(\mathcal{O})} & {H^1(L\otimes E_+)} & {H^1(L\otimes Q)} & 0 \\
        && {H^0(K)} \\
        && 0
        \arrow[from=1-3, to=2-3]
        \arrow[from=2-3, to=3-3]
        \arrow["{\delta'}", from=2-3, to=3-4]
        \arrow[from=3-1, to=3-2]
        \arrow[from=3-2, to=3-3]
        \arrow["r"', from=3-2, to=4-3]
        \arrow["\delta"', from=3-3, to=3-4]
        \arrow[from=3-3, to=4-3]
        \arrow["{\mu^*}", from=3-4, to=3-5]
        \arrow[from=3-5, to=3-6]
        \arrow[from=3-6, to=3-7]
        \arrow[from=4-3, to=5-3]
    \end{tikzcd}\]

    The horizontal exact sequence is \eqref{deformationLES} and the vertical short exact sequence comes from \eqref{torsionsequence} after tensoring by $L$. The map $r$ is the isomorphism induced by the sequence \eqref{d2ses}. This shows that the vertical sequence splits and we can write
    $$H^0(L\otimes Q)\simeq H^0(L|_{D_1})\oplus H^0(K).$$
    Also note that $H^1(L\otimes E_+)\simeq H^1(L\otimes Q)$ because they are both 1-dimensional, so $\mu^*=0$ and $\delta'$ is surjective. Under the identification $H^1(L\otimes Q)\simeq H^1(K)$ (see \eqref{torsionsequence}), the quadratic form $q$ is given by
    \begin{equation}\label{quadraticformgenus2}
        q(v+w)=\langle\delta'(w),v\rangle + \delta'(w)\cup w.
    \end{equation}
    Here, $v\in H^0(K)$ and $w\in H^0(L|_{D_1})$, and $\langle-,-\rangle:H^1(\O)\otimes H^0(K)\ra H^1(K)$ is the Serre duality pairing. Non-degeneracy of \eqref{quadraticformgenus2} is equivalent to $\{x:B(x,-)=0\}=0$, where $B$ is the associated bilinear form. Suppose $B(v+w,-)=0$. Then for every $v'\in H^0(K)$ we have $0=B(v+w,v')=\langle\delta'(w),v'\rangle$. The Serre duality pairing is non-degenerate and $\delta'$ is an isomorphism, so we get $w=0$. Now, for every $w'\in H^0(L|_{D_1})$ we have $0=B(v,w')=\langle\delta'(w'),v\rangle$, so $v=0$ for the same reasons. This shows that $q$ is non-degenerate, so $(L,\alpha)$ is an ordinary double point by Proposition \ref{ordinarydoublepoint}. 
\end{proof}

\begin{remark}
    One should expect the exceptional fibre of $\N_C^2(E_+)\ra J_C^2$ to contain an ordinary double point, because of the following local model: 
    $$\{(z_1,z_2,[x_0:x_1:x_2])\in\C^2\times\P^2 : z_1x_1+z_2x_2=0\}\ra \C^2_{z_1,z_2}.$$
    The fibre is $\P^1$ over $z=(z_1,z_2)\neq 0$, but over $z=0$ the fibre is $\P^2$. It is clear that $(0,0,[1:0:0])$ is an ordinary double point. Another way to think about this space is as the cone over the blowup $\Bl_0(\C^2)$.
\end{remark}

\begin{prop}\label{dequals2monodromy}
    Let $\E\ra C\times\DD$ be one of the two families $\E_{\pm}$ in Proposition \ref{semistablefamily}. Then the associated family of multi-vortex moduli spaces $\N_C^{2}(\E)\ra\DD$ is smooth. Its fibres over $t\neq 0$ are $\P^1$-bundles over $J_C^{2}$, the fibre over $t=0$ has one ordinary double point, and the monodromy map $\Upsilon$ is a Dehn twist.
\end{prop}

\begin{proof}
    The statements about the fibres follow from Theorem \ref{wobblystructure} and Lemma \ref{dequalstwo}. The smoothness of the total family follows a similar argument to Proposition \ref{W1smoothness}, since the non-zero nilpotent Higgs field is unique at a generic point in $\W_1$. Picard--Lefschetz theory then shows that the monodromy map is a Dehn twist.
\end{proof}

\begin{lemma}\label{nullhomologouslagrangians}
    Let $\pi:\N\ra X$ be a $\P^1$-bundle over a $2$-dimensional complex torus $X$. Then any 3-sphere embedded in $\N$ is nullhomologous. 
\end{lemma}

\begin{proof}
    Let $p:\widetilde{\N}\ra\N$ be the universal cover of $\N$. Then an embedding $i:S^3\hookrightarrow\N$ lifts to $\widetilde{\N}$ because $\pi_1(S^3)=0$. However, $\widetilde{\N}$ is a $\P^1$-bundle over $\C^2$, so $H_3(\widetilde{\N};\Z)=0$ and $i_*[S^3]=0$.
\end{proof}

\begin{lemma}\label{nofixedpoints}
    Let $\pi:\N\ra X$ be a $\P^1$-bundle over a $2$-dimensional complex torus $X$. Suppose $\Upsilon:\N\ra\N$ is an orientation-preserving diffeomorphism with $\Lambda(\Upsilon)=0$. Then $\Upsilon$ is homotopic to a homeomorphism with no fixed points. 
\end{lemma}

\begin{proof}
    We first show that $\Upsilon$ is homotopic to a a map covering a diffeomorphism of $X$. The Gysin sequence implies that $H^1(\N;\Z)\simeq H^1(X;\Z)$, so $\Upsilon^*$ can be viewed as a linear map $A:H^1(X;\Z)\ra H^1(X;\Z)$. Choose a diffeomorphism $f:X\ra X$ with $f^*=A$. Since
    $$[\N,X]=[\N,S^1]^4\simeq H^1(X;\Z)^4,$$
    the homotopy type of a map $\N\ra X$ is determined by its action on a basis of $H^1(X;\Z)$. This shows that $\pi\Upsilon\sim f\pi$ via a homotopy $H:\N\times I\ra X$. The fibration $\pi:\N\ra X$ satisfies the homotopy lifting property, so $H$ lifts to a homotopy $\tilde{H}:\N\times I\ra \N$ with $\tilde{H}_0=\Upsilon$ and $\pi\tilde{H}_1=\tilde{H}_1\pi$:
    \[\begin{tikzcd}
        {\N} & {\N} \\
        {\N\times I} & X
        \arrow["\Upsilon", from=1-1, to=1-2]
        \arrow[hook', from=1-1, to=2-1]
        \arrow["\pi", from=1-2, to=2-2]
        \arrow["{\tilde{H}}", dashed, from=2-1, to=1-2]
        \arrow["H"', from=2-1, to=2-2]
    \end{tikzcd}\]
    So, from now on we can assume that $\Upsilon$ covers $f:X\ra X$. As an abelian group, the cohomology of $\N$ decomposes as
    $$H^j(\N;\Z)=V_j\oplus\xi V_{j-2},$$
    where $V_j=\Lambda^jH^1(X;\Z)$ and $\xi\in H^2(\N;\Z)$ is the Chern class of the tautological line bundle. The induced map $\Upsilon^*$ acts on $V_j$ as $A_j=\Lambda^jA$, and there are $v\in V_2$ and $\lambda\in\Z$ such that $\Upsilon^*(\xi)=v+\lambda\xi$. As a matrix, this is
    $$\Upsilon^*=\begin{pmatrix}
        A_* & v\wedge A_{*-2}\\
        0 & \lambda A_{*-2}
    \end{pmatrix},$$
    which has trace $\tr(\Upsilon^*)=\tr(A_*) + \lambda\tr(A_{*-2})$. The condition that the lefschetz number of $\Upsilon$ vanishes can now be made more explicit:
    \begin{align*}
        0 &= 1 - \tr(A_1) + \tr(A_2) + \lambda - \tr(A_3) - \lambda\tr(A_1)\\
        &\hspace{4mm} + \tr(A_4) + \lambda\tr(A_2) - \lambda\tr(A_3) + \lambda\tr(A_4)\\
        &= (1+\lambda)\det(I-A).
    \end{align*}
    There are two cases to consider. Firstly, we may have $\det(I-A)=0$, which states that the Lefschetz number of $f:X\ra X$ vanishes. The main theorem of \cite{torifixedpoints} shows that any such map between tori is homotopic to a map without fixed points, implying that $\Upsilon$ itself has no fixed points. On the other hand, suppose that $\lambda=-1$ and $\det(I-A)\neq 0$. Then $f$ will have fixed points $x\in X$, and any fixed points of $\Upsilon$ must lie in $\N_x\simeq\P^1$. However since $\lambda=-1$, $\Upsilon|_{\N_x}$ is homotopic to the antipodal map. Then over a neighbourhood of $x$, $\Upsilon$ can be homotoped so that it restricts to the antipodal map on $\N_x$, in which case it has no fixed points. 
\end{proof}

The corollary below proves the final case $(d=2)$ of Theorem \ref{genus2theorem}.

\begin{cor}
    Let $C$ be a smooth genus $2$ curve. Then for any (admissible) loop in $\MM_C^{ss}\setminus\W_1$, the monodromy map $\Upsilon$ on $\N_C^2$ is homotopic to a diffeomorphism with no fixed points. In particular, $\#\Fix_{\mathfrak{s}}(\Upsilon)=0$ for all spin$^c$ structures of the form $\mathfrak{s}=(2,\gamma)$. 
\end{cor}

\begin{proof}
    Proposition \ref{dequals2monodromy} shows that the monodromy map $\Upsilon$ associated to any admissible loop is homotopic to a composition of Dehn twists. These all act trivially on homology by Lemma \ref{nullhomologouslagrangians}, so the Lefschetz number of $\Upsilon$ is $\chi(\N_C^2)=\chi(\P^1)\chi(J_C^2)=0$. Then Lemma \ref{nofixedpoints} shows that $\Upsilon$ is homotopic to map without fixed points, so all of its fixed point sets $\Fix_{\mathfrak{s}}(\Upsilon)$ are empty.
\end{proof}

\bibliographystyle{alpha}
\bibliography{ref}

\end{document}